\documentclass[12pt,unicode]{amsart}

\usepackage{cmap}
\usepackage[cp1251]{inputenc}
\usepackage[TS1,T2A]{fontenc}
\usepackage[russian]{babel}

\usepackage{geometry}
\usepackage{amsmath}
\usepackage{amsfonts}
\usepackage{amssymb}

\usepackage{pgf,tikz,circuitikz}
\usetikzlibrary{arrows}
\usepackage{graphicx}

\usepackage[pdftex,%
            colorlinks,linkcolor=red,citecolor=blue,urlcolor=blue,%
            pdfstartview=FitH,%
            pdfauthor={Sergey Dorichenko, Maxim Prasolov, Mikhail Skopenkov},%
            pdftitle={Dissections of a metal rectangle},%
            pdfkeywords={Dissection, Rectangle, Electrical network.}%
           ]{hyperref}

\newcommand{\mscomm}[1]{}
\newcommand{\No}{N}

\begin{document}

%\large

\title{Разрезания металлического прямоугольника}

\author{С. Дориченко, М. Прасолов и М. Скопенков}

\address{Школа 57 города Москвы}
\email{sdorichenko@gmail.com}

\address{Московский Государственный Университет}
\email{0x00002a@gmail.com}

\address{Институт проблем передачи информации Российской Академии Наук и Университет науки и техники Короля Абдулы}
\email{skopenkov@rambler.ru}

\date{}

\maketitle

\newcounter{problem}
\setcounter{problem}0

Задачи на разрезание наглядны и красивы, но иногда их совсем не просто решить.
%Задачи на разрезание с давних пор привлекали внимание благодаря своей наглядности и красоте, сочетающейся со сложностью решения. %\cite{D,FLR,FR,KeKi,LS,SuDi,Szegedy}.
С давних пор они вдохновляли дизайнеров и архитекторов. Ученые обратили на них внимание, когда обнаружилась их неожиданная связь с физикой и теорией вероятностей.
Об одной из таких задач и пойдет речь в этой статье.
%Интерес к этим задачам постоянно растет. Вот одна из них.
%Одна из них формулируется так:
%Интерес к ним растет в связи с обнаружением взаимосвязи этих задач с физикой, гармоническим анализом, теорией вероятностей.
%\cite{BeSch,BSST,CFP,DS,Fomin,Lovasz}.

%В этой статье мы ответим на вопрос, какие прямоугольники можно разрезать на квадраты. Ответ основан на физической интерпретации, использующей электрические цепи.
%Для ответа на первый вопрос нам потребуются электрические цепи постоянного тока, а на второй --- переменного.
%Для понимания статьи достаточно знакомства со школьной программой по математике и физике.

%\medskip

%{\bf ЧАСТЬ I. Разрезания прямоугольника на квадраты.}

%\smallskip

%Начнём с примера.

%Приведём пример разрезания прямоугольника на квадраты.

\medskip

\textbf{Какие прямоугольники можно разрезать на квадраты}
\smallskip

Прямоугольник размером $a\times b$, где $a$  и $b$  --- целые числа, легко разрезается на $a\cdot b$ одинаковых квадратов; см. рисунок~\ref{ristriv}.
%Разумеется, то же самое верно и для любого подобного ему прямоугольника.
Так же легко разрезать на равные квадраты прямоугольник с рациональным отношением сторон.

\begin{figure}[hb]
\definecolor{zzttqq}{rgb}{0.6,0.2,0}
\begin{tikzpicture}[line cap=round,line join=round,>=triangle 45,x=0.5cm,y=0.5cm]
\clip(-3.34,1.6) rectangle (3.4,5.96);
\fill[color=zzttqq,fill=zzttqq,fill opacity=0.1] (-3,5) -- (2,5) -- (2,2) -- (-3,2) -- cycle;
\draw [color=zzttqq] (-3,5)-- (2,5);
\draw [color=zzttqq] (2,5)-- (2,2);
\draw [color=zzttqq] (2,2)-- (-3,2);
\draw [color=zzttqq] (-3,2)-- (-3,5);
\draw [color=zzttqq] (-2,5)-- (-2,2);
\draw [color=zzttqq] (-1,5)-- (-1,2);
\draw [color=zzttqq] (0,5)-- (0,2);
\draw [color=zzttqq] (1,5)-- (1,2);
\draw [color=zzttqq] (-3,4)-- (2,4);
\draw [color=zzttqq] (-3,3)-- (2,3);
\draw[color=zzttqq] (-0.5,5.56) node {$a$};
\draw[color=zzttqq] (2.56,3.62) node {$b$};
\end{tikzpicture}
\caption{Прямоугольник $a\times b$ разрезается на $a\cdot b$ одинаковых квадратов.}
\label{ristriv}
\end{figure}

%Верно и обратное утверждение: если прямоугольник разрезан на одинаковые квадраты, то отношение его перпендикулярных сторон рационально.

Естественный вопрос: какие прямоугольники можно разрезать
на квадраты не обязательно одинакового размера.
%Поставим теперь следующий вопрос: каким может быть отношение сторон прямоугольника, разрезанного на \emph{не обязательно} одинаковые квадраты, которые не все одинаковы,
Оказывается, ответ тот же самый:

%Оказывается, это утверждение остается верным, если допустить в разрезании квадраты различных размеров:

\smallskip

\noindent{\bf Теорема Дена о разрезании прямоугольника.} {\it Если прямоугольник можно разрезать
на квадраты (не обязательно равные), то отношение длин его сторон рационально.}

%Оказывается, прямоугольники, которые можно разрезать на не обязательно
%одинаковые квадраты, тоже должны иметь рациональное отношение сторон.
%Оказывается, что этим примером исчерпываются все прямоугольники,
%разрезаемые на квадраты.

\smallskip

Эту теорему открыл Макс Ден в 1903 году.
Его доказательство было довольно сложным. Впоследствии появились более простые.
%\footnote{См., например, книгу И.М.~Яглома ``Как разрезать квадрат?'' М.: Наука, 1968. или обзор на сайте \url{http://www.squaring.net}.}. %\cite{BSST, Ya}.
Мы приведем одно из них, принадлежащее Р.Л.~Бруксу, К.А.Б.~Смиту, А.Г.~Стоуну и У.Т.~Татту\footnote{Увлекательный рассказ о том, как эти авторы придумали свой метод, можно прочитать в главе ``Квадрирование квадрата'' книги М.~Гарднера ``Математические головоломки и развлечения'', М.: Мир, 1999.}. Они придумали его, еще будучи студентами. Это доказательство основано на физической интерпретации, использующей электрические цепи. При этом физические соображения служат отправной точкой, а
само доказательство чисто математическое.
%проводится на математическом уровне строгости.
%Мы докажем такое более общее утверждение:

%Мы собираемся доказать теорему Дена в следующей более общей формулировке:

%\smallskip

%\noindent{\bf Обобщенная теорема Дена.}

%\noindent

%\noindent{\bf Теорема о разрезании прямоугольника.}
%{\it Пусть большой прямоугольник разрезан на меньшие прямоугольники.
%Тогда отношение сторон большого прямоугольника можно выразить  через
%отношения сторон меньших прямоугольников, пользуясь
%только операциями сложения, вычитания, умножения и деления}\footnote{В
%действительности, можно обойтись и без операции вычитания, но
%доказывать этого мы не будем.}.

%\smallskip

%Конечно, теорема Дена выводится из теоремы о разрезании прямоугольника.
%Действительно, пусть прямоугольник разрезан на квадраты. Отношение
%сторон любого квадрата равно $1$. С помощью указанных операций из
%числа $1$ можно получить только рациональное
%число. Значит, отношение сторон прямоугольника рационально.

%Далее {\it теоремой Дена} мы будем называть приведенное более общее
%утверждение.

Итак, пусть прямоугольник разрезан на квадраты. Чтобы найти отношение его сторон, достаточно найти стороны этих квадратов с точностью до пропорциональности. Покажем на примере, как это можно сделать.

\medskip

\textbf{Как найти стороны квадратов}

\smallskip

%Покажем на примере, как найти стороны квадратов, на которые разрезан прямоугольник.
%Расположим прямоугольник так, чтобы две его стороны были вертикальны, а две другие горизонтальны. Очевидно, что тогда стороны квадратов тоже вертикальны и горизонтальны.

%\mscomm{Начать с примера попроще?}

\begin{figure}[h]
\includegraphics[width=5cm]{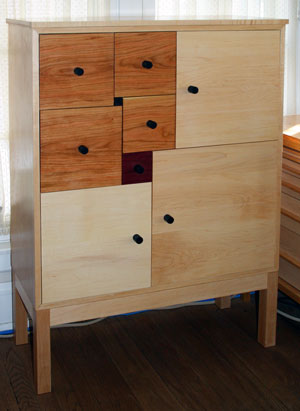}
\caption{Прямоугольный шкаф с квадратными полками.}
%Каково отношение ширины этого шкафа к его высоте (без учета ножек)?}
\label{ris1}
\end{figure}

%\noindent{\bf Пример 3.}
На рисунке~\ref{ris1} изображено
фото\footnote{Фото с сайта \url{http://www.mynl.com/ww/project11.html}} прямоугольного шкафа с квадратными полками. Представим себе, что мы хотим изготовить такой же шкаф. Для этого нам в первую очередь нужно узнать размеры полок. Просто измерить эти величины на фотографии не удастся, так как мы видим шкаф ``под углом'', а значит, истинные длины искажены.

\smallskip

%\bigskip

%\noindent{\it Решение задачи 1 (в конце номера).}
%Пусть прямоугольник
%разрезан на квадраты следующим образом (рис.~\ref{ris111}).
%Найдем отношение его сторон и покажем, что такая картинка действительно существует.

%Найдем отношение ширины шкафа к его высоте (без учета ножек).

%Поэтому мы воспользуемся более сложным методом.
Для того, чтобы найти эти размеры, занумеруем квадраты (полки) как показано на рисунке~\ref{ris111}.
Будем считать, что горизонтальная сторона прямоугольника (шкафа) равна $1$, а вертикальную сторону (без учета ножек) обозначим через $x$.
%Тогда искомое отношение сторон равно $R=1/x$.
Сторону квадрата $k$ обозначим через~$x_k$.

\begin{figure}[h]
  \definecolor{zzttqq}{rgb}{0.6,0.2,0}
\begin{tikzpicture}[line cap=round,line join=round,>=triangle 45,x=0.5cm,y=0.5cm]
\clip(-2.94,-4.13) rectangle (7.08,6.23);
\fill[color=zzttqq,fill=zzttqq,fill opacity=0.1] (6.9,-3.9) -- (6.9,6) -- (-2.7,6) -- (-2.7,-3.9) -- cycle;
\draw [color=zzttqq] (-2.7,-3.9)-- (6.9,-3.9);
\draw [color=zzttqq] (-2.7,6)-- (-2.7,-3.9);
\draw [color=zzttqq] (6.9,6)-- (-2.7,6);
\draw [color=zzttqq] (6.9,-3.9)-- (6.9,6);
\draw [color=zzttqq] (6.9,1.5)-- (1.5,1.5);
\draw [color=zzttqq] (1.5,1.5)-- (1.5,-3.9);
\draw [color=zzttqq] (2.4,6)-- (2.4,1.5);
\draw [color=zzttqq] (-2.7,0.3)-- (1.5,0.3);
\draw [color=zzttqq] (0.3,0.3)-- (0.3,3.6);
\draw [color=zzttqq] (1.5,1.5)-- (0.3,1.5);
\draw [color=zzttqq] (2.4,3.6)-- (0,3.6);
\draw [color=zzttqq] (0,6)-- (0,3.3);
\draw [color=zzttqq] (0.3,3.3)-- (-2.7,3.3);
\draw [color=zzttqq] (6.9,-3.9)-- (6.9,6);
\draw [color=zzttqq] (6.9,6)-- (-2.7,6);
\draw [color=zzttqq] (-2.7,6)-- (-2.7,-3.9);
\draw [color=zzttqq] (-2.7,-3.9)-- (6.9,-3.9);
\draw (0.08,3.71) node[anchor=north west] {1};
\draw (-1.64,2.05) node[anchor=north west] {2};
\draw (-1.74,5.15) node[anchor=north west] {3};
\draw (0.82,5.27) node[anchor=north west] {4};
\draw (1.02,2.99) node[anchor=north west] {5};
\draw (0.48,1.43) node[anchor=north west] {6};
\draw (3.88,-0.75) node[anchor=north west] {7};
\draw (-1.1,-1.33) node[anchor=north west] {8};
\draw (4.28,4.31) node[anchor=north west] {9};
\draw (0.18,3.43)-- (0.4,3.07);
\end{tikzpicture}
  \caption{Нумерация квадратов.%Разрезание прямоугольника на $9$ неравных квадратов
  }\label{ris111}
\end{figure}

К левой стороне прямоугольника
примыкают квадраты $2$, $3$ и $8$, откуда $x=x_2+x_3+x_8$.
К правой стороне квадрата $3$ примыкают квадраты $1$ и $4$:
$x_3=x_1+x_4$. Аналогично, $x_6+x_8=x_7$, $x_1+x_2=x_5+x_6$, $x_4+x_5=x_9$.
Условие для правой стороны прямоугольника мы не записываем, поскольку оно следует из предыдущих (получается сложением всех выписанных равенств). Сформулируем наше наблюдение, см.~рисунок~\ref{ris11111}:

\smallskip

\noindent{\bf Условие вертикальной стыковки.} {\it Для каждого вертикального разреза сумма сторон квадратов, примыкающих к разрезу слева, равна сумме сторон квадратов, примыкающих справа. Вертикальная сторона прямоугольника равна сумме сторон примыкающих к ней квадратов\footnote{
%Нетрудно проверить, что условие стыковки, записанное для правой стороны прямоугольника, является непосредственным следствием остальных условий вертикальной стыковки. Поэтому в дальнейшем мы выбросим его из нашей системы уравнений.\\
Для разрезаний, у которых в некоторых точках сходится сразу $4$ квадрата (как на рисунках~\ref{ristriv} или~\ref{riskontrprimer}), надо уточнить понятие разреза.
Покрасим все горизонтальные стороны квадратов, не лежащие на периметре прямоугольника, в зеленый цвет. Они объединятся в несколько зеленых отрезков, которые мы и назовем \emph{горизонтальными разрезами}.
%\emph{Горизонтальным разрезом} договоримся считать горизонтальный отрезок, проходящий по сторонам квадратов и ограниченный сторонами квадратов или периметром прямоугольника.
Вертикальные стороны квадратов, не лежащие на периметре, покрасим в оранжевый цвет. Полученные оранжевые отрезки делятся горизонтальными разрезами на части, именно эти части мы и назовем \emph{вертикальными разрезами}; см.~рисунок~\ref{riskontrprimer}.
%А \emph{вертикальным разрезом} будем считать вертикальный отрезок, проходящий по сторонам квадратов и ограниченный двумя горизонтальными разрезами или периметром прямоугольника; см.~рисунок~\ref{riskontrprimer}.
%Таким образом, по нашему соглашению, в разрезании на рисунке~\ref{riskontrprimer} один горизонтальный и два вертикальных разреза.
}.}
\smallskip

\begin{figure}[h]
\definecolor{ffwwqq}{rgb}{1,0.4,0}
\definecolor{zzttqq}{rgb}{0.6,0.2,0}
\begin{tikzpicture}[line cap=round,line join=round,>=triangle 45,x=0.5cm,y=0.5cm]
\clip(-3.05,0.19) rectangle (2.7,4.64);
\fill[color=zzttqq,fill=zzttqq,fill opacity=0.1] (0.3,3.3) -- (0.3,3.6) -- (0,3.6) -- (0,3.3) -- cycle;
\fill[color=zzttqq,fill=zzttqq,fill opacity=0.1] (0.3,1.5) -- (2.4,1.5) -- (2.4,3.6) -- (0.3,3.6) -- cycle;
\fill[color=zzttqq,fill=zzttqq,fill opacity=0.1] (1.5,0.3) -- (1.5,1.5) -- (0.3,1.5) -- (0.3,0.3) -- cycle;
\fill[color=zzttqq,fill=zzttqq,fill opacity=0.1] (0.3,0.3) -- (0.3,3.3) -- (-2.7,3.3) -- (-2.7,0.3) -- cycle;
\draw [color=zzttqq] (-2.7,0.3)-- (1.5,0.3);
\draw [line width=1.6pt,color=ffwwqq] (0.3,0.3)-- (0.3,3.6);
\draw [color=zzttqq] (1.5,1.5)-- (0.3,1.5);
\draw [color=zzttqq] (2.4,3.6)-- (0,3.6);
\draw [color=zzttqq] (0.3,3.3)-- (-2.7,3.3);
\draw [color=zzttqq] (0.3,3.6)-- (0,3.6);
\draw [color=zzttqq] (0,3.6)-- (0,3.3);
\draw [color=zzttqq] (0,3.3)-- (0.3,3.3);
\draw [color=zzttqq] (0.3,1.5)-- (2.4,1.5);
\draw [color=zzttqq] (2.4,1.5)-- (2.4,3.6);
\draw [color=zzttqq] (2.4,3.6)-- (0.3,3.6);
\draw [color=zzttqq] (1.5,0.3)-- (1.5,1.5);
\draw [color=zzttqq] (1.5,1.5)-- (0.3,1.5);
\draw [color=zzttqq] (0.3,0.3)-- (1.5,0.3);
\draw [color=zzttqq] (0.3,3.3)-- (-2.7,3.3);
\draw [color=zzttqq] (-2.7,3.3)-- (-2.7,0.3);
\draw [color=zzttqq] (-2.7,0.3)-- (0.3,0.3);
\draw (0.11,3.44)-- (-0.3,3.63);
\draw[color=zzttqq] (-0.69,3.8) node {$x_1$};
\draw[color=zzttqq] (1.52,2.5) node {$x_5$};
\draw[color=zzttqq] (0.96,0.82) node {$x_6$};
\draw[color=zzttqq] (-1.13,1.79) node {$x_2$};
\end{tikzpicture}
\caption{Условие вертикальной стыковки: $x_1+x_2=x_5+x_6$.}\label{ris11111}
\end{figure}

Заменяя слово ``вертикальный'' на ``горизонтальный'', а ``слева'' и ``справа'' --- на ``сверху'' и ``снизу'', мы получаем \textit{условие горизонтальной стыковки}, см.~рисунок~\ref{ris11112}.

%Глядя на схему разрезания, можно написать и уравнения, описывающие примыкание квадратов друг к другу горизонтальными сторонами:

%\smallskip
%\noindent{\bf Условие горизонтальной стыковки.} {\it Для каждого горизонтального разреза сумма горизонтальных сторон прямоугольников, примыкающих к разрезу сверху, равна сумме горизонтальных сторон прямоугольников, примыкающих снизу. Горизонтальная сторона большого прямоугольника равна сумме горизонтальных сторон примыкающих к ней прямоугольников.}
%\smallskip

\begin{figure}[h]
  \definecolor{qqccqq}{rgb}{0,0.8,0}
\definecolor{zzttqq}{rgb}{0.6,0.2,0}
\begin{tikzpicture}[line cap=round,line join=round,>=triangle 45,x=0.5cm,y=0.5cm]
\clip(0.02,-4.04) rectangle (7.12,6.18);
\fill[color=zzttqq,fill=zzttqq,fill opacity=0.1] (0.3,1.5) -- (0.3,0.3) -- (1.5,0.3) -- (1.5,1.5) -- cycle;
\fill[color=zzttqq,fill=zzttqq,fill opacity=0.1] (1.5,1.5) -- (1.5,-3.9) -- (6.9,-3.9) -- (6.9,1.5) -- cycle;
\fill[color=zzttqq,fill=zzttqq,fill opacity=0.1] (2.4,3.6) -- (0.3,3.6) -- (0.3,1.5) -- (2.4,1.5) -- cycle;
\fill[color=zzttqq,fill=zzttqq,fill opacity=0.1] (6.9,1.5) -- (6.9,6) -- (2.4,6) -- (2.4,1.5) -- cycle;
\draw [color=zzttqq] (6.9,-3.9)-- (6.9,6);
\draw [color=zzttqq] (1.5,1.5)-- (1.5,-3.9);
\draw [color=zzttqq] (2.4,6)-- (2.4,1.5);
\draw [color=zzttqq] (0.3,0.3)-- (0.3,3.6);
\draw [color=zzttqq] (0.3,1.5)-- (0.3,0.3);
\draw [color=zzttqq] (0.3,0.3)-- (1.5,0.3);
\draw [color=zzttqq] (1.5,0.3)-- (1.5,1.5);
\draw [color=zzttqq] (1.5,1.5)-- (1.5,-3.9);
\draw [color=zzttqq] (1.5,-3.9)-- (6.9,-3.9);
\draw [color=zzttqq] (6.9,-3.9)-- (6.9,1.5);
\draw [color=zzttqq] (2.4,3.6)-- (0.3,3.6);
\draw [color=zzttqq] (0.3,3.6)-- (0.3,1.5);
\draw [color=zzttqq] (2.4,1.5)-- (2.4,3.6);
\draw [color=zzttqq] (6.9,1.5)-- (6.9,6);
\draw [color=zzttqq] (6.9,6)-- (2.4,6);
\draw [color=zzttqq] (2.4,6)-- (2.4,1.5);
\draw [line width=1.6pt,color=qqccqq] (0.3,1.5)-- (6.9,1.5);
\draw[color=zzttqq] (1.06,0.9) node {$x_6$};
\draw[color=zzttqq] (4.44,-1.14) node {$x_7$};
\draw[color=zzttqq] (1.44,2.6) node {$x_5$};
\draw[color=zzttqq] (4.84,3.88) node {$x_9$};
\end{tikzpicture}
  \caption{Условие горизонтальной стыковки: $x_5+x_9=x_6+x_7$.} \label{ris11112}
\end{figure}

Из этого условия в нашем примере со шкафом получаем:  $1=x_3+x_4+x_9$, $x_4=x_1+x_5$, $x_1+x_3=x_2$, $x_5+x_9=x_6+x_7$, $x_2+x_6=x_8$. Условие для нижней стороны прямоугольника мы не записываем, поскольку оно следует из остальных.

Итак, осталось решить систему уравнений:
%все нахождение отношения сторон прямоугольника сводится к решению системы уравнений, в нашем случае:
\begin{equation*}\label{eq-styk}
\begin{matrix}
x=x_2+x_3+x_8, &x_3=x_1+x_4, &x_6+x_8=x_7, &x_1+x_2=x_5+x_6, &x_4+x_5=x_9,\\
x_3+x_4+x_9=1, &x_4=x_1+x_5, &x_1+x_3=x_2, &x_5+x_9=x_6+x_7, &x_2+x_6=x_8.
\end{matrix}
\end{equation*}

Такие уравнения называются \emph{линейными}.

\begin{figure}[h]
  \definecolor{qqccqq}{rgb}{0,0.8,0}
\definecolor{ffwwqq}{rgb}{1,0.4,0}
\definecolor{zzttqq}{rgb}{0.6,0.2,0}
\begin{tikzpicture}[line cap=round,line join=round,>=triangle 45,x=0.4cm,y=0.4cm]
\clip(-3.2,-0.3) rectangle (3.32,6.2);
\fill[color=zzttqq,fill=zzttqq,fill opacity=0.1] (-3,6) -- (3,6) -- (3,0) -- (-3,0) -- cycle;
\draw [color=zzttqq] (-3,6)-- (3,6);
\draw [color=zzttqq] (3,6)-- (3,0);
\draw [color=zzttqq] (3,0)-- (-3,0);
\draw [color=zzttqq] (-3,0)-- (-3,6);
\draw [line width=2.4pt,color=ffwwqq] (0,6)-- (0,0);
\draw [line width=2.4pt,color=qqccqq] (-3,3)-- (3,3);
\end{tikzpicture}
  \caption{В таком разрезании один горизонтальный и два вертикальных разреза.} \label{riskontrprimer}
\end{figure}

%Это --- система линейных уравнений относительно неизвестных
%$I,I_1,\dots,I_5$. Решив ее, мы выразим $I$ через $R_1,\dots,R_5$.
%Подставив это выражение в уравнение $RI=1$, найдем $R$.

%\smallskip

%\medskip
%{\bf Электрические цепи.}
%\smallskip

%Рассмотрим теперь более сложный пример: разрезание прямоугольника на
%попарно различные квадраты.

%В общем случае нужно действовать так же.

%Тем самым нахождение отношения сторон большого прямоугольника сводится к решению системы линейных уравнений.

%Посмотрим на условия стыковки с такой точки зрения. Горизонтальную сторону большого прямоугольника и отношения сторон меньших прямоугольников будем %считать данными. Представим горизонтальные стороны меньших прямоугольников как произведения вертикальной стороны на отношение сторон, тогда условия %стыковки --- это уравнения на неизвестные вертикальные стороны. Оказывается, что верна

%\smallskip

%\noindent{\bf Теорема единственности.} {
%\it Система уравнений, задаваемая условиями стыковки, в которой вертикальные стороны прямоугольников --- неизвестные,
%горизонтальная сторона большого прямоугольника равна $1$, а
%отношения сторон прямоугольников --- %и горизонтальная сторона большого прямоугольника
%данные числа, имеет единственное решение.}

%\smallskip

%Эту теорему мы докажем позже с помощью физической интерпретации. А пока будем двигаться дальше к доказательству теоремы о разрезании.

\medskip

{\bf Как решать систему линейных уравнений
\footnote{О решении линейных уравнений не раз рассказывалось в ``Кванте'', например, в статье В.~Гутенмахера в номере 1 за 1984 год.}}

\smallskip

Будем последовательно выражать неизвестные. В первом уравнении неизвестная $x$ выражена через другие неизвестные. Больше $x$ нигде не участвует, поэтому переходим ко второму уравнению. В нем неизвестная $x_3$ выражена через $x_1$ и $x_4$. Подставим это выражение в другие уравнения системы, содержащие неизвестную $x_3$ --- в первое, шестое и восьмое. Получим систему:
$$
\begin{matrix}
x=x_2+\mathbf{x_1+x_4}+x_8, &\mathbf{x_3=x_1+x_4}, &x_6+x_8=x_7, &x_1+x_2=x_5+x_6, &x_4+x_5=x_9,\\
\mathbf{x_1+2x_4}+x_9=1, &x_4=x_1+x_5, &\mathbf{2x_1+x_4}=x_2, &x_5+x_9=x_6+x_7, &x_2+x_6=x_8.
\end{matrix}
$$
Она равносильна исходной. Но теперь неизвестная $x_3$ участвует только во втором уравнении.

Перейдем к третьему уравнению. Подставляя выражение $x_7=x_6+x_8$ в девятое уравнение, получим систему, содержащую $x_7$ только в третьем уравнении:
$$
\begin{matrix}
x=x_2+x_1+x_4+x_8, &x_3=x_1+x_4, &\mathbf{x_6+x_8=x_7}, &x_1+x_2=x_5+x_6, &x_4+x_5=x_9,\\
x_1+2x_4+x_9=1, &x_4=x_1+x_5, &2x_1+x_4=x_2, &x_5+x_9=\mathbf{2x_6+x_8}, &x_2+x_6=x_8.
\end{matrix}
$$
Будем продолжать таким же образом дальше.
%Каждый раз будем выражать одну из неизвестных в очередном уравнении через остальные.
%Будем подставлять полученное выражение в другие уравнения. В полученной системе выбранная неизвестная будет присутствовать только в одном уравнении. После этого процесс продолжается\footnote{Вообще говоря, возможно, что после нашей подстановки все неизвестные в некотором уравнении сократятся, и оно примет вид $0=0$. В этом случае выбросим его из системы. Если же в результате подстановки получилось уравнение, в котором правая часть --- ненулевое число, а все коэффициенты в левой части нулевые, то это означает, что исходная система не имеет решений. %Иными словами, картинке вроде изображенной на рис.~\ref{ris111} не соответствует никакое разрезание прямоугольника на квадраты.
%}.

%На предпоследнем шаге мы получаем: \mscomm{???}
%$$
%\begin{matrix}
%x=33/32, &x_3=9/32, &x_7=9/16, &x_1=1/32, &x_4=1/4,\\
%x_9=15/32, &x_5=7/32, &x_2=5/16, &x_8=7/16, &x_6=1/8.
%\end{matrix}
%$$
В итоге мы получим систему ``уравнений'':
$$
\begin{matrix}
x=33/32, &x_3=9/32, &x_7=9/16, &x_1=1/32, &x_4=1/4,\\
x_9=15/32, &x_5=7/32, &x_2=5/16, &x_8=7/16, &x_6=1/8.
\end{matrix}
$$

Решение исходной системы найдено! Значения неизвестных $x_1,\dots,x_9$ --- это и есть стороны квадратов.
%, а отношение сторон прямоугольника равно $1/x=32/33$.
В нашем примере прямоугольник оказался разрезан на попарно различные квадраты.

\smallskip

\small

\refstepcounter{problem}
\noindent{\bf Задача \arabic{problem}.}\label{plane}
Докажите, что плоскость можно замостить попарно
различными квадратами, длины сторон которых ---
(а) рациональные; (б) целые числа.

\normalsize

\smallskip

\begin{figure}[h]
\includegraphics[width=5cm]{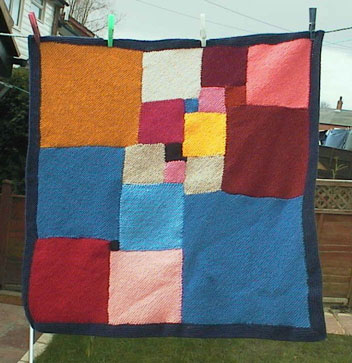}
\caption{Квадратное одеяло, сотканное из квадратных лоскутков.}
%\bf Нужно куда-то приплести эту картинку}
\label{ris11113}
\end{figure}

А можно ли \emph{квадрат} разрезать на попарно различные квадраты? Задача эта появилась в начале прошлого века, и оказалась очень сложной. Решили ее только спустя несколько десятилетий уже известные нам четыре студента и независимо от них Р.~Шпраг. Но если Р.~Шпраг использовал сложный перебор, то нашим студентам найти решение помогла физическая интерпретация. Потом было найдено много разных примеров, пример с наименьшим количеством квадратов изображен на рисунке~\ref{ris11113}.
%Оно иллюстрирует решение сложной старинной задачи: можно ли разрезать квадрат на различные квадраты?

\smallskip

\small

\noindent{\bf Задача \refstepcounter{problem}\arabic{problem}$^*$. \label{cube}}
Можно ли куб разрезать на несколько попарно-различных кубиков?

\normalsize

%Рассмотрим первое уравнение. Выразим одну из неизвестных в этом уравнении через остальные. Подставим полученное выражение в остальные уравнения %системы. В полученной системе выбранная неизвестная присутствует только в первом уравнении. (Возможно, что после нашей подстановки все неизвестные в %некотором уравнении сократились, и оно приняло вид $0=0$. В этом случае выбросим его из системы.)

%Рассмотрим второе уравнение. Выразим одну из неизвестных в этом уравнении через остальные. Подставим полученное выражение в остальные уравнения %системы. В полученной системе вторая выбранная нами неизвестная присутствует только во втором уравнении.

%Будем действовать таким образом, пока не закончатся уравнения. Если в итоге в каждом уравнении осталась ровно одна неизвестная, то решение находится %очевидным образом.

\medskip

%{\bf Если решение единственно, то оно рационально.}

{\bf Когда наш метод работает}

\smallskip

Итак, для шкафа мы нашли все интересующие нас размеры. Но будет ли так и для любого другого разрезания? Ясно, что если решение системы, построенной по условиям стыковки, единственно, то мы найдем его нашим методом. И, конечно же, оно будет рациональным: ведь коэффициенты системы рациональны, а мы, выражая неизвестные, использовали только сложение, вычитание, умножение и деление. Это простое наблюдение мы назовем так:

\smallskip
\noindent{\bf Теорема о решении системы.}%\footnote{Подробное доказательство этой теоремы можно найти в ...}
{\it Пусть система линейных уравнений с рациональными коэффициентами имеет единственное решение.
Тогда это решение состоит из рациональных чисел.}
%Тогда в этом решении каждая неизвестная выражается через коэффициенты уравнений с помощью четырёх арифметических действий.
\smallskip

Бывают системы линейных уравнений, у которых решение не единственно. Например, система
$$
x_1+x_2=0,\quad
x_1+x_3=1;\qquad
$$
имеет бесконечно много решений: в качестве $x_1$ можно взять любое число $t$, в качестве $x_2$ число $-t$, а в качестве $x_3$ число $1-t$. У нее есть и иррациональные решения (когда $t$ иррационально).

%Например, алгоритм Гаусса--Жордана приводит систему уравнений
%$$
%x_1+x_2=0,\quad
%x_1+x_3=1;\qquad
%\text{ к виду }\qquad
%x_1=-x_2,\quad
%x_3=1+x_2.
%$$
%При этом неизвестную $x_2$ мы ``не успели'' выразить через остальные. Для каждого значения $t$ этой переменной тройка $x_1=-t, x_2=t, x_3=1+t$ является решением системы.

%Мы доказали следующую теорему:

%Заметим, что коэффициента в системе уравнений, которую мы составили, равны либо $\pm1$, либо отношениям сторон меньших прямоугольников.

Совершенно не очевидно, что условий стыковки достаточно, чтобы найти стороны квадратов, то есть, что система, построенная по реальному разрезанию, имеет единственное решение.
%Однако это так: если система построена по реальному разрезанию, то решение всегда единственно.
Оказывается, что это всегда так: наш метод позволяет \emph{однозначно} восстановить все размеры по фотографии разрезания (если мы считаем горизонтальную сторону прямоугольника равной $1$). %(разумеется, с точностью до одновременного умножения всех размеров на одно и то же число).
%прямоугольника, разрезанного на квадраты (например, для рисунка~\ref{ris11113}), наш метод позволяет \emph{однозначно} восстановить все истинные размеры, то есть решение системы всегда единственно.
Мы докажем это с помощью физической интерпретации. А теорема Дена о разрезании прямоугольника отсюда сразу следует по теореме о решении системы.
%Итак, для доказательства теоремы Дена о разрезании прямоугольника остается показать, что система уравнений, построенная по условиям стыковки, имеет единственное решение.
%Это мы сделаем с помощью физической интерпретации. %При этом все необходимые понятия из физики мы напомним.

%Доказательство этого факта мы пока отложим, а сейчас перейдем к его следствиям. \mscomm{Непонятно, к каким следствиям: структура неясна.}

\smallskip

\small

\refstepcounter{problem}
\noindent{\bf Задача \arabic{problem}. \label{rooms}}
Архитектор нарисовал план квартиры. На плане (рисунок~\ref{ris9})
показано, как должны примыкать комнаты друг к другу, но их размеры искажены. Можно ли сделать все комнаты квадратными?
%разбить прямоугольник на квадраты так, как показано на рисунке~\ref{ris9}.

\normalsize

\smallskip

\begin{figure}[h]
    \definecolor{zzttqq}{rgb}{0.6,0.2,0}
\begin{tikzpicture}[line cap=round,line join=round,>=triangle 45,x=0.3cm,y=0.3cm]
\clip(-4.3,-1.3) rectangle (3.44,6.3);
\fill[color=zzttqq,fill=zzttqq,fill opacity=0.1] (-4,6) -- (3,6) -- (3,-1) -- (-4,-1) -- cycle;
\draw [color=zzttqq] (-4,6)-- (3,6);
\draw [color=zzttqq] (3,6)-- (3,-1);
\draw [color=zzttqq] (3,-1)-- (-4,-1);
\draw [color=zzttqq] (-4,-1)-- (-4,6);
\draw [color=zzttqq] (-1,6)-- (-1,2);
\draw [color=zzttqq] (-4,2)-- (0,2);
\draw [color=zzttqq] (-1,3)-- (3,3);
\draw [color=zzttqq] (0,3)-- (0,-1);
\end{tikzpicture}
    \caption{План квартиры.} \label{ris9}
\end{figure}

%\begin{figure}[h]
%    \caption{Искаженное и истинное изображения прямоугольника, разрезанного на квадраты}\label{ris11114}
%\end{figure}

%\smallskip
%\small
%\noindent{\bf Задача \refstepcounter{problem}\arabic{problem}.} Докажите, что если система линейных уравнений имеет единственное решение, то в этом решении каждая неизвестная выражается через коэффициенты уравнений с помощью четырёх арифметических действий.

\normalsize

\bigskip
\bigskip
\bigskip
\bigskip

{\bf Физическая интерпретация.}

\smallskip

Оказывается, каждому разрезанию прямоугольника на квадраты можно сопоставить электрическую цепь. %см~рисунок~\ref{rismain}.
Если мы найдем токи в этой электрической цепи, то мы найдем и стороны квадратов. Но обо всем по порядку.

%Мы объяснили часть названия статьи --- в рассматриваемых вопросах речь
%идет о {\it разрезаниях прямоугольника}.
%Но одно слово из названия осталось без объяснения --- почему именно
%''{\it металлического}''?

%Мы объяснили совпадение выражений в рассмотренных примерах с формулами для сопротивления простейших цепей.
%Однако, наше обоснование было, что называется, \textit{на физическом уровне строгости}. Теперь мы собираемся дать \textit{математически строгое} обоснование.
%Построим {\it математическую модель} идеальной электрической цепи.
% и покажем, как выглядят разрезания прямоугольника с точки зрения электрических цепей. Так, например, мы увидим, что единственность системы уравнений, построенной по правилам стыковки, означает, что электрический ток в цепи распределяется единственным образом.
%Напомним, что \emph{граф} --- это некоторое множество точек, называемых \emph{вершинами}, соединенных между собой линиями, называемыми \emph{ребрами}.
Мы будем
рассматривать \emph{математическую модель} электрической цепи\footnote{Желающим подробно разобраться
в физике происходящего рекомендуем статью ``Правила Кирхгофа'' в Кванте \No1 за 1985 год.}.
%\mscomm{Попросить у физиков ссылку на хороший учебник}}.
Вместо физических законов и опытных фактов у нас будут определения, аксиомы, теоремы.

С математической точки зрения \textit{электрическая цепь} --- это связный плоский граф, каждому ребру которого сопоставлено некоторое положительное число, причем концы одного из ребер отмечены знаками ``$+$'' и ``$ -  $''.
%Будем использовать следующие понятия.
Ребро с отмеченными концами называется {\it батарейкой}, остальные --- {\it резисторами}. Число, сопоставленное батарейке, называется \emph{напряжением} батарейки, а числа, сопоставленные резисторам, --- их \emph{сопротивлениями}. Вершины графа называются \emph{узлами}, отмеченные узлы батарейки --- \emph{положительной} и \emph{отрицательной} \emph{клеммами}.

%При желании читатель может воспринимать приводимую физическую интерпретацию просто как
%удобный язык, и не вникать в физический смысл происходящего.
%Везде в дальнейшем без ущерба для понимания можно считать, что {\it цепь} --- это некоторый граф, нарисованный на плоскости, а токи и сопротивления %--- это некоторые числа, стоящие на рёбрах этого графа.
%подразумевать разрезание, под {\it резистором}~---
%прямоугольник, под {\it сопротивлением}~--- отношение его сторон, под {\it мощностью}~--- площадь
%, под {\it силой тока}~--- длины вертикальных сторон
%и т.д\footnote{Строго говоря, доказательство сформулированной ниже теоремы единственности всё же требует перехода на язык цепей%.}.
%Достоинство физической интерпретации в том, что на этом языке все шаги решения
%возникают естественным образом, в то время
%как на языке исходной задачи решение выглядело бы %необъяснимым
%набором трюков.
По разрезанию цепь строится так;
%Покажем, как сопоставить разрезанию прямоугольника на квадраты электрическую цепь; %(как математический объект);
см.~рисунок~\ref{rismain}.
%на рисунках~\ref{ris3} и~\ref{rismain} показаны примеры.

%\mscomm{Дать сначала рисунок без стрелок и чисел}

\begin{figure}[h]
\begin{tabular}{cc}
\definecolor{ffqqqq}{rgb}{1,0,0}
\definecolor{qqwwff}{rgb}{0,0.4,1}
\definecolor{ffwwqq}{rgb}{1,0.4,0}
\definecolor{zzttqq}{rgb}{0.6,0.2,0}
\begin{tikzpicture}[line cap=round,line join=round,>=triangle 45,x=0.5cm,y=0.5cm]
\clip(-4.04,-4.1) rectangle (8.2,10.02);
\fill[color=zzttqq,fill=zzttqq,fill opacity=0.1] (6.9,-3.9) -- (6.9,6) -- (-2.7,6) -- (-2.7,-3.9) -- cycle;
\draw [color=zzttqq] (-2.7,-3.9)-- (6.9,-3.9);
\draw [line width=1.6pt,color=ffwwqq] (-2.7,6)-- (-2.7,-3.9);
\draw [color=zzttqq] (6.9,6)-- (-2.7,6);
\draw [line width=1.6pt,color=ffwwqq] (6.9,-3.9)-- (6.9,6);
\draw [color=zzttqq] (6.9,1.5)-- (1.5,1.5);
\draw [line width=1.6pt,color=ffwwqq] (1.5,1.5)-- (1.5,-3.9);
\draw [line width=1.6pt,color=ffwwqq] (2.4,6)-- (2.4,1.5);
\draw [color=zzttqq] (-2.7,0.3)-- (1.5,0.3);
\draw [line width=1.6pt,color=ffwwqq] (0.3,0.3)-- (0.3,3.6);
\draw [color=zzttqq] (1.5,1.5)-- (0.3,1.5);
\draw [color=zzttqq] (2.4,3.6)-- (0,3.6);
\draw [line width=1.6pt,color=ffwwqq] (0,6)-- (0,3.3);
\draw [color=zzttqq] (0.3,3.3)-- (-2.7,3.3);
\draw [color=qqwwff] (6.9,1.08)-- (1.5,-1.2);
\draw [color=qqwwff] (6.9,1.08)-- (2.4,3.76);
\draw [color=qqwwff] (1.5,-1.2)-- (-2.7,1.04);
\draw [color=qqwwff] (0.3,1.8)-- (2.4,3.76);
\draw [color=qqwwff] (1.5,-1.2)-- (0.3,1.8);
\draw [color=qqwwff] (0,4.66)-- (2.4,3.76);
\draw [color=qqwwff] (0,4.66)-- (-2.7,1.04);
\draw [color=qqwwff] (0.3,1.8)-- (-2.7,1.04);
\draw [color=qqwwff] (0,4.66)-- (0.3,1.8);
\draw [shift={(2.09,4.14)},color=qqwwff]  plot[domain=-0.57:3.72,variable=\t]({1*5.71*cos(\t r)+0*5.71*sin(\t r)},{0*5.71*cos(\t r)+1*5.71*sin(\t r)});
\fill [color=ffqqqq] (6.9,1.08) circle (2.0pt);
\draw[color=ffqqqq] (7.34,1.18) node {$-$};
\fill [color=ffqqqq] (1.5,-1.2) circle (2.0pt);
\fill [color=ffqqqq] (2.4,3.76) circle (2.0pt);
\fill [color=ffqqqq] (-2.7,1.04) circle (2.0pt);
\draw[color=ffqqqq] (-3.22,1.08) node {$+$};
\fill [color=ffqqqq] (0,4.66) circle (2.0pt);
\fill [color=ffqqqq] (0.3,1.8) circle (2.0pt);
\end{tikzpicture} &
\ctikzset{bipoles/length=.8cm}
\begin{circuitikz}[scale=0.6, american]
\draw[color=brown]
(-2.7,1) to[generic, *-*] (0,4.66)
      to[generic, *-*] (2.4,3.76)
      to[generic, *-*] (6.9,1)
(-2.7,1) to[generic, *-*] (1.5,-1.2)
      to[generic, *-*] (6.9,1)
(-2.7,1) to[generic, *-*] (0.3,1.8)
      to[generic, *-*] (2.4,3.76)
(0,4.66) to[generic, *-*]  (0.3,1.8)
      to[generic, *-*] (1.5,-1.2)
(6.9,1)  to[short, *-] (6.9,6.0) 
      to[battery, -] (-2.7,6.0) 
      to[short, -*] (-2.7,1)  
%(-2.7,1)  to[short, *-] (-2.7,6.0)     
%      to[battery, -, v^<=$ v_1$] (6.9,6.0) 
%      to[short, -*] (6.9,1)      
;\end{circuitikz}
\end{tabular}
\caption{Построение электрической цепи по разрезанию (слева). Общепринятое изображение электрической цепи (справа).
%Резисторы принято изображать прямоугольниками, а батарейку --- специальным значком (справа).
}\label{rismain}
\end{figure}

На каждой вертикальной линии разреза отметим по точке --- это будут узлы будущей электрической цепи. На вертикальных сторонах прямоугольника выберем по клемме. Их мы отметим знаками ``$+$'' и ``$-$'' и соединим с батарейкой (``$+$'' на левой стороне, ``$-$'' на правой).

Каждый квадрат ограничен слева и справа двумя вертикальными разрезами. В электрической цепи его изображением служит резистор, соединяющий два узла на этих разрезах (узлы могут оказаться на продолжениях сторон квадрата).
%Одна из этих клемм изображает левую сторону маленького прямоугольника, а другая --- правую.

Сопротивление каждого резистора положим равным\footnote{Мы раз и навсегда фиксируем систему единиц: сопротивления будем измерять в килоомах, напряжения в вольтах, токи --- в миллиамперах. В дальнейшем единицы измерения не указываются.}
 $1$.
Напряжение батарейки также положим равным $1$.

%Клеммы, соответствующие вертикальным сторонам большого прямоугольника, соединяются с разными полюсами батарейки (на рисунке~\ref{rismain} сама батарейка не изображена).
%Нам удобно будет положить напряжение батарейки  численно равным длине горизонтальной стороны большого прямоугольника.
Нужная нам электрическая цепь построена.

\medskip

{\bf Как найти токи в электрической цепи}
%\footnote{На эту тему см.~также статьи А.~Варламова ``Правила Кирхгофа'', О.~Ляшко ``Почему не уменьшится сопротивление'' в Кванте \No1 (1985) и циклы задач Летних конференциях Турнира городов: М.~Прасолов, М.~Скопенков, Б.~Френкин, ``Инварианты многоугольников'', \url{http://www.turgor.ru/lktg/2007/1/index.php}, Д.~Баранов, М.~Скопенков, А.~Устинов, ``Случайные блуждания и электрические цепи'', \url{http://olympiads.mccme.ru/lktg/2010/4/index.htm}.}.

\smallskip
Теперь объясним, что такое \emph{токи} в электрической цепи и как их можно найти.

%Вернёмся к нашей математической модели.
Занумеруем резисторы, как показано на рисунке~\ref{rismain2} (то есть так же, как соответствующие квадраты).
%Нарисуем на ребре с отмеченными концами стрелку справа налево, а на остальных ребрах --- слева направо
Нарисуем на каждом резисторе стрелку слева направо,
а на батарейке --- справа налево, то есть, от отрицательного полюса к положительному\footnote{Мы нарисовали предполагаемые направления тока.
%Если мы ошиб\"емся с выбором направления тока на каком-то ребре, то в итоге мы получим верную величину тока с обратным знаком.
Читателя может смутить, что в одном из ребер ток направлен от ``минуса'' к ``плюсу''.
%Это может смутить читателя, знающего, что ``ток в цепи идет от плюса к минусу''.
Но это действительно так:
ток в {\it резисторах} идет от ``плюса'' к ``минусу'', а вот в {\it батарейке} --- наоборот.
}.
%припишем току в каждом резисторе одно из двух направлений (мы можем ошибиться с выбором направления, но тогда у нас в итоге для силы тока получится то же значение со знаком минус). Для определенности в рассматриваемом примере направим ток на резисторах слева направо.
%Наоборот, на батарейке нарисуем стрелку справа налево\footnote{Это может смутить читателя, знающего, что ``ток в цепи идет от плюса к минусу''. Ток в {\it цепи} действительно идет в этом направлении, а вот в самой {\it батарейке} --- наоборот}.
%В цепи появился электрический ток. Как было сказано, ток в резисторе определяется направлением и величиной ({\it силой тока}). %(такие объекты называют {\it векторами}).
Электрическая цепь делит плоскость на части. % (одна из которых бесконечна).
Обходя границу любой %конечной
части по часовой стрелке, получим замкнутую цепочку ребер, называемую {\it контуром}\footnote{
Для простоты будем считать, что контур не проходит ни через какое ребро дважды. Это не всегда так, см.~рисунок~\ref{risgran}. Однако в дальнейшем мы увидим, что это так для любой цепи, построенной по разрезанию прямоугольника.}.

\begin{figure}[h]
\definecolor{ffqqqq}{rgb}{1,0,0}
\definecolor{qqwwff}{rgb}{0,0.4,1}
\begin{tikzpicture}[line cap=round,line join=round,>=triangle 45,x=0.5cm,y=0.5cm]
\clip(-4.46,-1.58) rectangle (8.26,9.74);
\draw [color=qqwwff] (6.9,1.08)-- (1.5,-1.2);
\draw [color=qqwwff] (6.9,1.08)-- (2.4,3.76);
\draw [color=qqwwff] (1.5,-1.2)-- (-2.7,1.04);
\draw [color=qqwwff] (0.3,1.8)-- (2.4,3.76);
\draw [color=qqwwff] (1.5,-1.2)-- (0.3,1.8);
\draw [color=qqwwff] (0,4.66)-- (2.4,3.76);
\draw [color=qqwwff] (0,4.66)-- (-2.7,1.04);
\draw [color=qqwwff] (0.3,1.8)-- (-2.7,1.04);
\draw [color=qqwwff] (0,4.66)-- (0.3,1.8);
\draw [->,color=qqwwff] (-2.7,1.04) -- (0,4.66);
\draw [->,color=qqwwff] (0,4.66) -- (2.4,3.76);
\draw [->,color=qqwwff] (2.4,3.76) -- (6.9,1.08);
\draw [->,color=qqwwff] (0,4.66) -- (0.3,1.8);
\draw [->,color=qqwwff] (-2.7,1.04) -- (0.3,1.8);
\draw [->,color=qqwwff] (0.3,1.8) -- (2.4,3.76);
\draw [->,color=qqwwff] (-2.7,1.04) -- (1.5,-1.2);
\draw [->,color=qqwwff] (0.3,1.8) -- (1.5,-1.2);
\draw [->,color=qqwwff] (1.5,-1.2) -- (6.9,1.08);
\draw [shift={(2.09,3.92)},color=qqwwff]  plot[domain=-0.53:3.68,variable=\t]({1*5.59*cos(\t r)+0*5.59*sin(\t r)},{0*5.59*cos(\t r)+1*5.59*sin(\t r)});
\draw [->,color=qqwwff] (-2.89,1.39) -- (-2.7,1.04);
\draw[color=qqwwff] (4.28,-0.6) node {7};
\fill [color=ffqqqq] (6.9,1.08) circle (2.0pt);
\draw[color=ffqqqq] (7.34,1.18) node {$-$};
\draw[color=qqwwff] (5,2.98) node {9};
\draw[color=qqwwff] (-0.88,-0.86) node {8};
\draw[color=qqwwff] (1.7,2.48) node {5};
\fill [color=ffqqqq] (1.5,-1.2) circle (2.0pt);
\draw[color=qqwwff] (1.28,0.8) node {6};
\fill [color=ffqqqq] (2.4,3.76) circle (2.0pt);
\draw[color=qqwwff] (1.38,4.84) node {4};
\draw[color=qqwwff] (-1.32,4.04) node {3};
\fill [color=ffqqqq] (-2.7,1.04) circle (2.0pt);
\draw[color=ffqqqq] (-3.22,1.08) node {$+$};
\draw[color=qqwwff] (-1.12,1.94) node {2};
\fill [color=ffqqqq] (0,4.66) circle (2.0pt);
\fill [color=ffqqqq] (0.3,1.8) circle (2.0pt);
\draw[color=qqwwff] (0.72,3.38) node {1};
\end{tikzpicture}
\caption{Нумерация резисторов и выбор направлений на резисторах и батарейке.}
\label{rismain2}
\end{figure}

\begin{figure}[h]
\definecolor{qqwwff}{rgb}{0,0.4,1}
\definecolor{ffqqqq}{rgb}{1,0,0}
\begin{tikzpicture}[line cap=round,line join=round,>=triangle 45,x=0.3cm,y=0.3cm]
\clip(-3.3,-3.87) rectangle (8.74,6.41);
\draw [color=qqwwff] (0.88,-3.49)-- (-2.7,5.57);
\draw [color=qqwwff] (-2.7,5.57)-- (7.81,5.48);
\draw [color=qqwwff] (7.81,5.48)-- (0.88,-3.49);
\draw [line width=1.6pt,color=qqwwff] (0.88,-3.49)-- (3.22,1.85);
\draw [color=qqwwff] (3.22,1.85)-- (4.32,3.66);
\draw [color=qqwwff] (4.32,3.66)-- (1.28,3.59);
\draw [color=qqwwff] (3.22,1.85)-- (1.28,3.59);
\draw [->,dash pattern=on 2pt off 2pt] (-1.82,5.06) -- (6.93,5);
\draw [->,dash pattern=on 2pt off 2pt] (6.93,5) -- (2.25,-1.14);
\draw [->,dash pattern=on 2pt off 2pt] (2.25,-1.14) -- (3.66,1.63);
\draw [->,dash pattern=on 2pt off 2pt] (3.66,1.63) -- (5.68,4.14);
\draw [->,dash pattern=on 2pt off 2pt] (5.68,4.14) -- (-0.11,4.01);
\draw [->,dash pattern=on 2pt off 2pt] (-0.11,4.01) -- (2.62,1.7);
\draw [->,dash pattern=on 2pt off 2pt] (2.62,1.7) -- (1.06,-2.04);
\draw [->,dash pattern=on 2pt off 2pt] (1.06,-2.04) -- (-1.82,5.06);
\fill [color=ffqqqq] (0.88,-3.49) circle (1.5pt);
\fill [color=ffqqqq] (-2.7,5.57) circle (1.5pt);
\draw[color=ffqqqq] (-2.72,4.14) node {$+$};
\fill [color=ffqqqq] (7.81,5.48) circle (1.5pt);
\draw[color=ffqqqq] (7.75,4.07) node {$-$};
\fill [color=ffqqqq] (3.22,1.85) circle (1.5pt);
\fill [color=ffqqqq] (4.32,3.66) circle (1.5pt);
\fill [color=ffqqqq] (1.28,3.59) circle (1.5pt);
\end{tikzpicture}
\caption{Контур, проходящий по ребру дважды.}
\label{risgran}
\end{figure}

%Мы рассматриваем лишь такие электрические цепи, которые можно нарисовать на плоскости так, чтобы внутренности ребер не пересекались; хотя рассматриваемый метод переносится с небольшими изменениями на произвольные электрические цепи.

%Немного поговорим о физике происходящего. Напомним, что сопротивление резистора --- это характеристика самого резистора, не зависящая от цепи, в которую он подключен. Мы считаем батарейку идеальной, то есть обладающей нулевым внутренним сопротивлением. В этом случае напряжение батарейки также не зависит от подключаемой цепи.

%Когда цепь, содержащая батарейку, замыкается, в ней начинает течь электрический ток.
%Ток в каждом резисторе и батарейке характеризуется направлением и величиной ({\it силой тока}). Сила тока может быть и нулевой. Покажем, как найти силы тока в электрической цепи.
%Суммарная сила тока, протекающего через цепь, пропорциональна напряжению батарейки. Данное соотношение было открыто экспериментально и известно как \emph{закон Ома}. Величина, обратная коэффициенту пропорциональности, называется \emph{сопротивлением} цепи. Ток измеряется в Амперах, напряжение --- в Вольтах, а сопротивление --- в Омах.
\emph{Сила тока} через $k$-й резистор  --- это просто некоторое действительное число $I_k$, сопоставленное резистору. %(оно меняет знак при смене направления стрелки на резисторе на противоположное).
{\it Сила тока} через батарейку --- это некоторое действительное число $I$. {\it Напряжение} на резисторе --- это произведение силы тока на его сопротивление. (А для батарейки напряжение вообще от тока не зависит. Такая батарейка в физике называется \emph{идеальной}.)
%{\it Сопротивление} цепи --- это отношение напряжения батарейки к силе тока через нее.
Силы тока определяются следующими аксиомами (правилами), проиллюстрированными на~рисунках~\ref{rispravilo1}, \ref{rispravilo2}, \ref{rispravilo2a}:
%Обозначим через $I_1$, $I_2$, \dots, $I_9$ --- неизвестные силы тока в резисторах, а через $I$ --- неизвестную силу тока в батарейке. Составим уравнения на числа $I_1, I_2, \dots, I_9$, пользуясь следующими правилами.
%\footnote{Правила Кирхгофа подтверждаются экспериментально}

%Физическое наблюдение состоит в том, что силы тока в резисторах удовлетворяют уравнениям Кирхгофа:

\smallskip

\noindent{\bf Первое правило Кирхгофа.}
{\it В каждом узле сумма входящих токов равна сумме выходящих.}
%\footnote{ Нетрудно проверить, что первое правило Кирхгофа для самой правой клеммы электрической цепи (отмеченной знаком ``$-$'') является непосредственным следствием первого правила Кирхгофа для всех остальных клемм. Поэтому в дальнейшем мы выбрасываем из нашей системы уравнение, полученное из первого правила Кирхгофа для этой клеммы.}
%(рис.\ref{ris6} слева).}

\begin{figure}[h]
\definecolor{ffqqqq}{rgb}{1,0,0}
\definecolor{qqwwff}{rgb}{0,0.4,1}
\begin{tikzpicture}[line cap=round,line join=round,>=triangle 45,x=0.5cm,y=0.5cm]
\clip(-2.98,-1.52) rectangle (2.68,4.9);
\draw [color=qqwwff] (0.3,1.8)-- (2.4,3.76);
\draw [color=qqwwff] (1.5,-1.2)-- (0.3,1.8);
\draw [color=qqwwff] (0.3,1.8)-- (-2.7,1.04);
\draw [color=qqwwff] (0,4.66)-- (0.3,1.8);
\draw [->,color=qqwwff] (0,4.66) -- (0.3,1.8);
\draw [->,color=qqwwff] (-2.7,1.04) -- (0.3,1.8);
\draw [->,color=qqwwff] (0.3,1.8) -- (2.4,3.76);
\draw [->,color=qqwwff] (0.3,1.8) -- (1.5,-1.2);
\draw[color=qqwwff] (1.78,2.44) node {$I_5$};
\fill [color=ffqqqq] (1.5,-1.2) circle (2.0pt);
\draw[color=qqwwff] (1.48,0.68) node {$I_6$};
\fill [color=ffqqqq] (2.4,3.76) circle (2.0pt);
\fill [color=ffqqqq] (-2.7,1.04) circle (2.0pt);
\draw[color=qqwwff] (-1.3,2) node {$I_2$};
\fill [color=ffqqqq] (0,4.66) circle (2.0pt);
\fill [color=ffqqqq] (0.3,1.8) circle (2.0pt);
\draw[color=qqwwff] (-0.5,3.48) node {$I_1$};
\end{tikzpicture}
\caption{Первое правило Кирхгофа: $I_1+I_2=I_5+I_6$.}
\label{rispravilo1}
\end{figure}

\smallskip

Для нашего примера получаем такие уравнения:
$I=I_2+I_3+I_8, I_3=I_1+I_4, I_6+I_8=I_7, I_1+I_2=I_5+I_6, I_4+I_5=I_9.$
%$I=I_1+I_4$; $I_1=I_2+I_5$; $I_3=I_2+I_4$.
(Мы не записываем уравнение для самой правой клеммы, поскольку оно непосредственно следует из остальных.)

%\smallskip

%Назовём {\it контуром} цепочку резисторов, получающуюся при обходе против часовой стрелки границы любой части, на которые цепь делит плоскость. %(рис.\ref{ris6} справа).

\smallskip

\noindent{\bf Второе правило Кирхгофа.} {\it Для любого контура сумма напряжений на резисторах (с  соответствующими знаками) равна напряжению батарейки (с соответствующими знаком),
если контур содержит батарейку, а иначе равна нулю. Напряжение на резисторе берется со знаком ``$+$'', если направление стрелки на резисторе %или батарейке
совпадает с направлением обхода контура, а иначе со знаком ``$ - $''.
%если эти направления противоположны.
Так же определяется знак для батарейки.}

\begin{figure}[h]
\definecolor{ffqqqq}{rgb}{1,0,0}
\definecolor{qqwwff}{rgb}{0,0.4,1}
\begin{tikzpicture}[line cap=round,line join=round,>=triangle 45,x=0.5cm,y=0.5cm]
\clip(-0.08,-1.4) rectangle (7.24,4.04);
\draw [color=qqwwff] (6.9,1.08)-- (1.5,-1.2);
\draw [color=qqwwff] (6.9,1.08)-- (2.4,3.76);
\draw [color=qqwwff] (0.3,1.8)-- (2.4,3.76);
\draw [color=qqwwff] (1.5,-1.2)-- (0.3,1.8);
\draw [->,color=qqwwff] (2.4,3.76) -- (6.9,1.08);
\draw [->,color=qqwwff] (0.3,1.8) -- (2.4,3.76);
\draw [->,color=qqwwff] (0.3,1.8) -- (1.5,-1.2);
\draw [->,color=qqwwff] (1.5,-1.2) -- (6.9,1.08);
\draw [dash pattern=on 2pt off 2pt] (3.08,1.32) circle (0.51cm);
\draw [->] (3.18,2.34) -- (3.44,2.28);
\draw[color=qqwwff] (4.52,-0.6) node {$I_7$};
\fill [color=ffqqqq] (6.9,1.08) circle (2.0pt);
\draw[color=qqwwff] (5.24,2.98) node {$I_9$};
\draw[color=qqwwff] (1.94,2.48) node {$I_5$};
\fill [color=ffqqqq] (1.5,-1.2) circle (2.0pt);
\draw[color=qqwwff] (1.52,0.8) node {$I_6$};
\fill [color=ffqqqq] (2.4,3.76) circle (2.0pt);
\fill [color=ffqqqq] (0.3,1.8) circle (2.0pt);
\end{tikzpicture}
\caption{Второе правило Кирхгофа
 для контура без батарейки:
%:
$I_5+I_9-I_6-I_7=0$ (учтено, что все сопротивления в нашем примере равны $1$)
.}
\label{rispravilo2}
\end{figure}

\begin{figure}[h]
\definecolor{qqwwff}{rgb}{0,0.4,1}
\definecolor{ffqqqq}{rgb}{1,0,0}
\begin{tikzpicture}[line cap=round,line join=round,>=triangle 45,x=0.5cm,y=0.5cm]
\clip(-4.46,0.3) rectangle (8.12,9.74);
\draw [color=qqwwff] (6.9,1.08)-- (2.4,3.76);
\draw [color=qqwwff] (0,4.66)-- (2.4,3.76);
\draw [color=qqwwff] (0,4.66)-- (-2.7,1.04);
\draw [->,color=qqwwff] (-2.7,1.04) -- (0,4.66);
\draw [->,color=qqwwff] (0,4.66) -- (2.4,3.76);
\draw [->,color=qqwwff] (2.4,3.76) -- (6.9,1.08);
\draw [shift={(2.09,3.92)},color=qqwwff]  plot[domain=-0.53:3.68,variable=\t]({1*5.59*cos(\t r)+0*5.59*sin(\t r)},{0*5.59*cos(\t r)+1*5.59*sin(\t r)});
\draw [->,color=qqwwff] (-2.89,1.39) -- (-2.7,1.04);
\draw [dash pattern=on 3pt off 3pt] (2.08,6.86) circle (0.64cm);
\draw [->] (2.06,8.14) -- (2.32,8.12);
\fill [color=ffqqqq] (6.9,1.08) circle (2.0pt);
\draw[color=ffqqqq] (7.34,1.18) node {$-$};
\draw[color=qqwwff] (5.24,2.98) node {$I_9$};
\fill [color=ffqqqq] (2.4,3.76) circle (2.0pt);
\draw[color=qqwwff] (1.62,4.84) node {$I_4$};
\draw[color=qqwwff] (-1.08,4.04) node {$I_3$};
\fill [color=ffqqqq] (-2.7,1.04) circle (2.0pt);
\draw[color=ffqqqq] (-3.22,1.08) node {$+$};
\fill [color=ffqqqq] (0,4.66) circle (2.0pt);
\end{tikzpicture}
\caption{Второе правило Кирхгофа для контура
с батарейкой:
%:
$-I_3-I_4-I_9=-1$ (учтено, что напряжение %
батарейки
в нашем примере равно $1$).}
\label{rispravilo2a}
\end{figure}

\smallskip

Поскольку у нас напряжение батарейки  равно $1$  и все сопротивления равны $1$, получаем уравнения:
$-I_3-I_4-I_9=-1, I_4-I_1-I_5=0, I_1+I_3-I_2=0, I_5+I_9-I_6-I_7=0, I_2+I_6-I_8=0.$
(Мы не записываем уравнение для контура вокруг всей цепи, поскольку оно непосредственно следует из остальных.)
Решая полученную систему уравнений, находим все силы токов:
$$
\begin{matrix}
I=33/32, &I_3=9/32, &I_7=9/16, &I_1=1/32, &I_4=1/4,\\
I_9=15/32, &I_5=7/32, &I_2=5/16, &I_8=7/16, &I_6=1/8.
\end{matrix}
$$
%Поэтому общее сопротивление электрической цепи равно
%$R=1/I=32/33$.

У нас есть аксиомы (правила Кирхгофа), которые мы заимствовали из физики, а все остальные утверждения об электрических цепях мы выводим из них чисто математически.

%Из школьного курса физики
%Известны следующие примеры вычисления общего сопротивления (сравните их с примерами~1 и~2). Мы предлагаем читателю самостоятельно вывести эти правила из наших аксиом.

%\smallskip

%\noindent{\bf Пример 1$'$.} Пусть два резистора $R_1$ и $R_2$ соединены {\it последовательно} (см.~нижнюю часть рисунка~\ref{ris3}а). Тогда сопротивление полученной цепи равно $R_1+R_2$.

%\smallskip

%\noindent{\bf Пример 2$'$.} Пусть два резистора $R_1$ и $R_2$ соединены {\it параллельно} (см.~нижнюю часть рисунка~\ref{ris3}б). Тогда сопротивление полученной цепи вычисляется по формуле $\frac{R_1R_2}{R_1+R_2}$.

%\smallskip

%\begin{figure}

  % Рисунка пока нет

 % \caption{Последовательное и параллельное соединения}\label{ris5}

%\end{figure}

%Не любая электрическая цепь ``сводится'' к последовательно и параллельно соединенным резисторам (например, см.~нижнюю часть рисунка~\ref{ris3}в). Покажем на примере рисунка~\ref{rismain}, как искать сопротивление ``плоской'' цепи в общем случае. Сопротивления всех резисторов в цепи считаем равным $1$ Ом, напряжение батарейки --- равным $1$ Вольт.

\smallskip

\small

%Следующая небольшая серия задач позволит нам установить следствие, которое пригодится в дальнейшем.

%\smallskip
%\noindent{\bf Определение.} Потенциалом клеммы называется величина, равная сумме напряжений (с такими же знаками, как во втором правиле Кирхгофа) на всех резисторах некоторого пути, соединяющего положительный полюс батарейки с данной клеммой.

%\noindent{\bf Задача \refstepcounter{problem}\arabic{problem}.} Докажите, что потенциал клеммы не зависит от выбора пути.

\noindent{\bf Задача \refstepcounter{problem}\arabic{problem}\label{closed chain}.}
Выведите из второго правила Кирхгофа более общее правило, которое получается, если заменить в формулировке контур на любую замкнутую цепочку ребер (не проходящую ни через какую вершину дважды).
%Докажите, что второе правило Кирхгофа остается справедливым, если в его формулировке заменить контур на любую замкнутую цепочку ребер (без повторений).
%%%%
%Рассмотрим в электрической цепи любой замкнутый путь, идущий по стрелкам (не обязательно контур). Выведите из второго правила Кирхгофа, что сумма напряжений на резисторах пути равна напряжению батарейки, если путь содержит батарейку, а иначе равна нулю.
%Выведите из второго правила Кирхгофа, что для любой замкнутой цепочки ребер (а не только для контура) сумма напряжений на резисторах (с подходящими знаками) равна напряжению батарейки (с подходящим знаком), если цепочка содержит батарейку, а иначе равна нулю.
%%%
%сумма напряжений (с подходящими знаками) на любой замкнутой цепочке резисторов (а не только контуре) равна нулю.
%Докажите, что если во втором правиле Кирхгофа ``любой контур'' заменить на ``любой замкнутый путь в графе'', то мы получим систему уравнений, эквивалентную исходной.

\normalsize

%\smallskip

\medskip

{\bf Правила Кирхгофа совпадают с условиями стыковки}

\smallskip

%Мы собираемся сопоставить разрезанию
%прямоугольника на прямоугольники электрическую цепь из резисторов с одним источником тока (рис.~\ref{ris2}).

%Мы сопоставили разрезанию прямоугольника электрическую цепь.

Удивительным образом, правила Кирхгофа дают нам ту же самую систему уравнений на силы токов, что и условия стыковки на длины сторон квадратов! Докажем это.
%Итак, покажем, что правила Кирхгофа для токов совпадают с условиями стыковки на длины сторон квадратов.
%Покажем,
%что система уравнений на длины сторон прямоугольников, построенная по условиям стыковки, совпадает с системой уравнений на силы тока, построенной по правилам Кирхгофа.
%%%%%%Положим силу тока на резисторе численно равной стороне соответствующего квадрата.
%Выберем направление тока в каждом резисторе от клеммы, соответствующей левой стороне прямоугольника, к клемме, сопоставленной правой стороне того же прямоугольника.
%Найдём решение системы уравнений Кирхгофа. Подадим напряжение, численно равное длине горизонтальной стороны большого прямоугольника. Предположим, %что ток в проводнике идёт в направлении от точки, соответствующей левой стороне прямоугольника, к точке, сопоставленной правой стороне того же %прямоугольника.

%(см. рисунок~\ref{ris22222}).

%\begin{figure}


%  \input 8.tex

 % \caption{Как по разрезанию построить цепь}\label{ris2}

%\end{figure}

Рассмотрим первое правило Кирхгофа. Зафиксируем вертикальный разрез и соответствующий ему узел.
Входящие в узел токи соответствуют сторонам квадратов, примыкающим к разрезу слева, а выходящие из узла --- сторонам квадратов справа, см.~рисунок~\ref{ris5555}. Значит, первое правило Кирхгофа в этом узле для токов совпадает с правилом вертикальной стыковки.

%Рассмотрим первое правило Кирхгофа. Зафиксируем вертикальный разрез и соответствующую ему клемму. Входящие в клемму токи соответствуют сторонам квадратов, примыкающим к разрезу слева, а выходящие из клеммы --- сторонам квадратов справа, см.~рисунок~\ref{ris5555}. Значит, первое правило Кирхгофа в этой клемме для токов совпадает с правилом вертикальной стыковки.

\begin{figure}
\definecolor{ffwwqq}{rgb}{1,0.4,0}
\definecolor{ffqqqq}{rgb}{1,0,0}
\definecolor{qqwwff}{rgb}{0,0.4,1}
\definecolor{zzttqq}{rgb}{0.6,0.2,0}
\begin{tikzpicture}[line cap=round,line join=round,>=triangle 45,x=0.5cm,y=0.5cm]
\clip(-2.94,-1.42) rectangle (2.72,5);
\fill[color=zzttqq,fill=zzttqq,fill opacity=0.1] (0.3,3.6) -- (2.4,3.6) -- (2.4,1.5) -- (0.3,1.5) -- cycle;
\fill[color=zzttqq,fill=zzttqq,fill opacity=0.1] (-2.7,3.3) -- (0.3,3.3) -- (0.3,0.3) -- (-2.7,0.3) -- cycle;
\fill[color=zzttqq,fill=zzttqq,fill opacity=0.1] (0,3.6) -- (0.3,3.6) -- (0.3,3.3) -- (0,3.3) -- cycle;
\fill[color=zzttqq,fill=zzttqq,fill opacity=0.1] (0.3,1.5) -- (1.5,1.5) -- (1.5,0.3) -- (0.3,0.3) -- cycle;
\draw [color=zzttqq] (2.4,6)-- (2.4,0.6);
\draw [color=zzttqq] (1.5,1.5)-- (1.5,-3.9);
\draw [color=zzttqq] (0.3,0.3)-- (0.3,3.6);
\draw [color=zzttqq] (0,6)-- (0,3.3);
\draw [color=zzttqq] (0.3,3.3)-- (-2.7,3.3);
\draw [color=qqwwff] (0.3,1.8)-- (2.4,3.76);
\draw [color=qqwwff] (1.5,-1.2)-- (0.3,1.8);
\draw [color=qqwwff] (0.3,1.8)-- (-2.7,1.04);
\draw [color=qqwwff] (0,4.66)-- (0.3,1.8);
\draw [->,color=qqwwff] (0,4.66) -- (0.3,1.8);
\draw [->,color=qqwwff] (-2.7,1.04) -- (0.3,1.8);
\draw [->,color=qqwwff] (0.3,1.8) -- (2.4,3.76);
\draw [->,color=qqwwff] (0.3,1.8) -- (1.5,-1.2);
\draw [color=zzttqq] (0.3,3.6)-- (2.4,3.6);
\draw [color=zzttqq] (2.4,3.6)-- (2.4,1.5);
\draw [color=zzttqq] (2.4,1.5)-- (0.3,1.5);
\draw [color=zzttqq] (0.3,1.5)-- (0.3,3.6);
\draw [color=zzttqq] (-2.7,3.3)-- (0.3,3.3);
\draw [color=zzttqq] (0.3,3.3)-- (0.3,0.3);
\draw [color=zzttqq] (0.3,0.3)-- (-2.7,0.3);
\draw [color=zzttqq] (-2.7,0.3)-- (-2.7,3.3);
\draw [color=zzttqq] (0,3.6)-- (0.3,3.6);
\draw [color=zzttqq] (0.3,3.6)-- (0.3,3.3);
\draw [color=zzttqq] (0.3,3.3)-- (0,3.3);
\draw [color=zzttqq] (0,3.3)-- (0,3.6);
\draw [color=zzttqq] (0.3,1.5)-- (1.5,1.5);
\draw [color=zzttqq] (1.5,1.5)-- (1.5,0.3);
\draw [color=zzttqq] (1.5,0.3)-- (0.3,0.3);
\draw [color=zzttqq] (0.3,0.3)-- (0.3,1.5);
\draw [line width=1.6pt,color=ffwwqq] (0.3,3.58)-- (0.3,0.3);
\draw[color=qqwwff] (1.78,2.44) node {$I_5$};
\fill [color=ffqqqq] (1.5,-1.2) circle (2.0pt);
\draw[color=qqwwff] (0.38,-0.58) node {$I_6$};
\fill [color=ffqqqq] (2.4,3.76) circle (2.0pt);
\fill [color=ffqqqq] (-2.7,1.04) circle (2.0pt);
\draw[color=qqwwff] (-1.3,2) node {$I_2$};
\fill [color=ffqqqq] (0,4.66) circle (2.0pt);
\fill [color=ffqqqq] (0.3,1.8) circle (2.0pt);
\draw[color=qqwwff] (-0.7,4.28) node {$I_1$};
\end{tikzpicture}
\caption{Первое правило Кирхгофа и условие вертикальной стыковки.}
\label{ris5555}
\end{figure}

Рассмотрим второе правило Кирхгофа. Возьмем любой горизонтальный разрез. Ясно, что резисторы, соответствующие примыкающим к нему квадратам, образуют контур, см.~рисунок~\ref{risserge}. Квадраты сверху образуют верхнюю часть контура, а квадраты снизу --- нижнюю. Поскольку все сопротивления единичны,
напряжение на каждом резисторе равно силе тока на нем.
Значит, правило горизонтальной стыковки для нашего разреза совпадает со вторым правилом Кирхгофа.

Наоборот, возьмем любой контур. Его самый левый узел
соответствует некоторому вертикальному разрезу. К нему примыкают справа два квадрата,
соответствующие двум выходящим из  узла резисторам контура. Рассмотрим горизонтальный разрез, к которому примыкают эти квадраты. Снова, все квадраты, примыкающие к разрезу сверху, образуют верхнюю часть нашего контура,
а все примыкающие снизу --- нижнюю. Значит, второе правило Кирхгофа для нашего контура совпадает с правилом горизонтальной стыковки.

%Рассмотрим второе правило Кирхгофа. Возьмем любой контур нашей цепи. Его самая левая клемма  соответствует некоторому вертикальному разрезу. К нему примыкают справа два квадрата, соответствующие двум выходящим из клеммы резисторам контура. Рассмотрим горизонтальный разрез, к которому примыкают эти квадраты. Верхняя часть нашего контура тогда соответствует квадратам, примыкающим к горизонтальному разрезу сверху, а нижняя --- квадратам, примыкающим снизу.

\begin{figure}
\definecolor{qqccqq}{rgb}{0,0.8,0}
\definecolor{ffqqqq}{rgb}{1,0,0}
\definecolor{qqwwff}{rgb}{0,0.4,1}
\definecolor{zzttqq}{rgb}{0.6,0.2,0}
\begin{tikzpicture}[line cap=round,line join=round,>=triangle 45,x=0.5cm,y=0.5cm]
\clip(-0.11,-4.39) rectangle (7.48,6.38);
\fill[color=zzttqq,fill=zzttqq,fill opacity=0.1] (0.3,0.3) -- (1.5,0.3) -- (1.5,1.5) -- (0.3,1.5) -- cycle;
\fill[color=zzttqq,fill=zzttqq,fill opacity=0.1] (1.5,-3.9) -- (6.9,-3.9) -- (6.9,1.5) -- (1.5,1.5) -- cycle;
\fill[color=zzttqq,fill=zzttqq,fill opacity=0.1] (0.3,3.6) -- (0.3,1.5) -- (2.4,1.5) -- (2.4,3.6) -- cycle;
\fill[color=zzttqq,fill=zzttqq,fill opacity=0.1] (2.4,1.5) -- (6.9,1.5) -- (6.9,6) -- (2.4,6) -- cycle;
\draw [color=zzttqq] (6.9,6)-- (6.9,-3.6);
\draw [color=zzttqq] (6.9,-3.9)-- (6.9,6);
\draw [color=zzttqq] (1.5,1.5)-- (1.5,-3.9);
\draw [color=zzttqq] (0.3,0.3)-- (0.3,3.6);
\draw [color=qqwwff] (6.9,1.08)-- (1.5,-1.2);
\draw [color=qqwwff] (6.9,1.08)-- (2.4,3.76);
\draw [color=qqwwff] (0.3,1.8)-- (2.4,3.76);
\draw [color=qqwwff] (1.5,-1.2)-- (0.3,1.8);
\draw [->,color=qqwwff] (2.4,3.76) -- (6.9,1.08);
\draw [->,color=qqwwff] (0.3,1.8) -- (2.4,3.76);
\draw [->,color=qqwwff] (0.3,1.8) -- (1.5,-1.2);
\draw [->,color=qqwwff] (1.5,-1.2) -- (6.9,1.08);
\draw [line width=1.6pt,color=qqccqq] (0.3,1.5)-- (6.9,1.5);
\draw [color=zzttqq] (0.3,0.3)-- (1.5,0.3);
\draw [color=zzttqq] (1.5,0.3)-- (1.5,1.5);
\draw [color=zzttqq] (0.3,1.5)-- (0.3,0.3);
\draw [color=zzttqq] (1.5,-3.9)-- (6.9,-3.9);
\draw [color=zzttqq] (6.9,-3.9)-- (6.9,1.5);
\draw [color=zzttqq] (1.5,1.5)-- (1.5,-3.9);
\draw [color=zzttqq] (0.3,3.6)-- (0.3,1.5);
\draw [color=zzttqq] (2.4,1.5)-- (2.4,3.6);
\draw [color=zzttqq] (2.4,3.6)-- (0.3,3.6);
\draw [color=zzttqq] (6.9,1.5)-- (6.9,6);
\draw [color=zzttqq] (6.9,6)-- (2.4,6);
\draw [color=zzttqq] (2.4,6)-- (2.4,1.5);
\draw[color=qqwwff] (4.69,-0.85) node {$I_7$};
\fill [color=ffqqqq] (6.9,1.08) circle (2.0pt);
\draw[color=qqwwff] (5.2,3.35) node {$I_9$};
\draw[color=qqwwff] (1.75,2.33) node {$I_5$};
\fill [color=ffqqqq] (1.5,-1.2) circle (2.0pt);
\draw[color=qqwwff] (0.55,-0.58) node {$I_6$};
\fill [color=ffqqqq] (2.4,3.76) circle (2.0pt);
\fill [color=ffqqqq] (0.3,1.8) circle (2.0pt);
\end{tikzpicture}
\caption{Второе правило Кирхгофа и условие горизонтальной стыковки.}
%\caption{Контуру соответствует горизонтальный разрез}
\label{risserge}
\end{figure}

%Поскольку все сопротивления единичны, напряжение на каждом резисторе равно силе тока на нем.

Итак, правила Кирхгофа совпадают с условиями стыковки. Значит ли это, что токи совпадают с длинами сторон квадратов? Да, но только если наша система уравнений
имеет \emph{только одно} решение.
%токи \emph{однозначно} определяются правилами Кирхгофа.
В этом случае, изготовив по разрезанию электрическую цепь, длины сторон квадратов можно было бы найти... просто измерив токи!

%Возьмем любой контур в нашей цепи. Пусть $m$ --- самый левый узел этой цепи, он соответствует правой стороне некоторого прямоугольника $M$,  и $n$ --- %самый правый (соответствует левой стороне прямоугольника $N$).
%Пройдем по верхней части контура. Соседним резисторам соответствуют соседние
%прямоугольники разбиения. Они примыкают друг к другу так, что их нижние горизонтальные
%стороны находятся на одной прямой. Ведь иначе нашлись бы два соседних прямоугольника
%$A$ и $B$, такие что у $A$ нижняя горизонтальная сторона выше, чем у $B$.
%Тогда к прямоугольникам $A$ и $B$ примыкает еще прямоугольник $C$ (см.рис.),
%и  соответствующий участок контура выглядит так:

%Получается, что внутрь нашей грани отходит ребро, что невозможно (докажите!).

%Тогда от левой стороны  прямоугольника  $M$ к правой стороне прямоугольника $N$ идет горизонтальный разрез. Верхняя часть нашего контура соответствует %прямоугольникам, примыкающим к этому разрезу сверху, а нижняя ---
%прямоугольникам, примыкающим к разрезу снизу. Второе правило Кирхгофа означает просто,
%что суммы длин горизонтальных сторон тех прямоугольников, которые примыкают к разрезу
%сверху, и тех, которые примыкают снизу, одинаковы --- это очевидно, поскольку обе суммы
%просто равны длине разреза.

%Значит, по теореме единственности, длины горизонтальных сторон прямоугольников разрезания численно равны силам токов на соответствующих резисторах. А %тогда сопротивление цепи численно равно отношению сторон большого прямоугольника.

\medskip
\textbf{Единственность распределения токов в электрической цепи}
\smallskip

%Теперь покажем, что силы тока (а значит, и стороны квадратов) однозначно находятся из нашей системы:

%\smallskip

\noindent{\bf Теорема единственности.} {\it
 Пусть сопротивления всех резисторов цепи положительны. Тогда
система уравнений, построенная по правилам Кирхгофа, в которой силы тока --- неизвестные, а напряжение батарейки и сопротивления резисторов известны, имеет не более одного решения.}

\smallskip

На ``физическом уровне строгости'' эта теорема почти очевидна. Пусть решений два. Вычтем одно из другого. Тогда напряжение батарейки станет нулевым, а ток не везде будет равен нулю, чего не бывает.
% Предположим, что наша система имеет два решения. Тогда их разность удовлетворяет той же системе, только с нулевым напряжением батарейки.
%Но при нулевом напряжении батарейки тока нет. Значит, разность решений равна нулю --- то есть, решения совпадают.

С точки зрения математики это объяснение нельзя считать доказательством. Нельзя исключить возможность, что наша система уравнений имеет какие-то ``посторонние'' решения, которые не реализуются в ``реальной'' электрической цепи.
%Даже если мы примем на веру, что в ``реальной'' электрической цепи при нулевом напряжении ток нулевой, остается возможность, что наша система уравнений имеет какие-то ``посторонние'' ненулевые решения.
%А нам нужно, чтобы именно {\it наша система} не имела ненулевых решений.
Да и в нашем рассуждении мы нигде не использовали, что сопротивления всех резисторов строго положительны. А без этого предположения теорема неверна: в кольце из сверхпроводника (то есть резистора с нулевым сопротивлением) может течь ненулевой ток при нулевом напряжении!

Вот как можно математически строго доказать
теорему единственности:

%Предположим, что при данных напряжении батарейки, сопротивлениях и выбранных направлениях резисторов система уравнений, построенная по правилам Кирхгофа, имеет два различных решения.

\smallskip

\noindent {\bf Доказательство теоремы единственности.}
% Пусть $I_1,I_2,\dots$ и $I_1',I_2',\dots$ --- два разных решения. Тогда их разность $I_1-I_1',I_2-I_2',\dots$
Пусть есть два решения. Первое будем обозначать $I_1$, $I_2$, $\ldots$, второе --- $J_1$, $J_2$, $\ldots$.
Наша цель --- доказать, что их разность  $I_1-J_1$, $I_2-J_2$, $\ldots$ нулевая.

Рассмотрим любое уравнение нашей системы. Пусть, например, оно записано для  узла, изображенного на рисунке~\ref{rispravilo1}. Подставив в него первое решение, получим: $I_1+I_2=I_5+I_6$. Подставив второе, получим: $J_1+J_2=J_5+J_6$. Вычтем одно равенство из другого: $(I_1-J_1)+(I_2-J_2)=(I_5-J_5)+(I_6-J_6)$. Получается, что разность наших решений удовлетворяет тому же самому уравнению. Так будет и для уравнения, записанного для  любого другого узла или любого контура, не содержащего батарейку.

Пусть теперь уравнение записано для контура, содержащего батарейку, например, для контура на рисунке~\ref{rispravilo2a}. Подставляя в это уравнение наши решения, получим равенства $I_3+I_4+I_9=1$ и $J_3+J_4+J_9=1$. Вычтем одно равенство из другого: $(I_3-J_3)+(I_4-J_4)+(I_9-J_9)=0$. Получается, что разность наших решений удовлетворяет тому же самому уравнению, только с нулевой правой частью. Но в правой части исходного уравнения стояло напряжение батарейки. Получаем, что разность наших решений подчиняется правилам Кирхгофа для той же цепи, только с нулевым напряжением батарейки.

Теорема единственности свелась к такому утверждению:

\smallskip
\noindent{\bf Принцип техники безопасности.} \emph{Если напряжение батарейки равно нулю, то и все силы тока в электрической цепи нулевые.}

\smallskip

\noindent{\bf Доказательство.}
%Докажем его.
%\noindent{\bf Доказательство теоремы единственности в случае, когда напряжение батарейки равно нулю.}
Пусть в цепи есть ненулевые токи. Если сила тока на каких-то ребрах отрицательна, то поменяем  на каждом из них направление стрелки, знак силы тока и напряжения. Ясно, что правила Кирхгофа по-прежнему будут выполняться, а все силы тока станут неотрицательны.
%Если одновременно изменить направление стрелки на некотором ребре и поменять знак силы тока на нем, то правила Кирхгофа по-прежнему будут выполняться.
%Сделаем это для всех ребер, где сила тока отрицательна.
Начнем движение с ребра,
на котором сила тока ненулевая,
и будем двигаться в направлении стрелок. Из первого правила Кирхгофа следует, что мы можем неограниченно продолжать движение (ведь если у вершины есть положительный входящий ток, то есть и выходящий). Рано или поздно мы впервые вернемся в вершину, в которой уже побывали.
Значит, мы получим замкнутую цепочку
ребер, на которых сила тока неотрицательна, причем хотя бы на одном из них она больше нуля.
% на каждом резисторе которого напряжение положительно.
По задаче~\ref{closed chain} получаем противоречие со вторым правилом Кирхгофа, потому что напряжение батарейки равно нулю. %Следовательно, система уравнений, построенная по правилам Кирхгофа, имеет единственное решение.
Принцип техники безопасности, а вместе с ним и теорема единственности, доказаны.
%Рассмотрим произвольный резистор с ненулевой силой тока. Будем двигаться по нему в выбранном направлении. Мы попадем в некоторую клемму цепи. В нее входит ненулевой ток по резистору, по которому мы в нее пришли. Сумма входящих токов равна сумме выходящих, значит, найдется другой резистор, по которому ненулевой ток выходит из данной клеммы. Продолжим движение по этому новому резистору в выбранном на нем направлении. Продолжая таким образом, в конце концов мы получим цикл, не содержащий батарейку, на каждом резисторе которого напряжение положительно. По задаче \ref{closed chain} получаем противоречие со вторым правилом Кирхгофа. Следовательно, система уравнений, построенная по правилам Кирхгофа, имеет единственное решение. Теорема доказана.

%\smallskip

%Читатель, конечно, заметил, что для рассматриваемого примера система уравнений, построенная по правилам Кирхгофа, идентична системе, построенной ранее по условиям стыковки. Это проявление общей закономерности, которую мы установим в следующем пункте.

%\normalsize

%Мы докажем этот факт позже.
%Если смотреть на эту систему уравнений формально, то условие положительности сопротивлений существенно.

%Заметим, что общее сопротивление цепи не зависит от напряжения батарейки: если подать на цепь другое
%напряжение, то все силы токов можно умножить на отношение нового напряжения батарейки к старому, и тогда для новых сил тока правила Кирхгофа будут %выполнены. Осталось заметить, что отношение напряжения батарейки к общей силе тока осталось прежним.

\smallskip

\small

\noindent{\bf Задача \refstepcounter{problem}\arabic{problem}$^*$.
\label{polozhitelnost}} Выведите из правил Кирхгофа, что если напряжение батарейки положительно, то сила тока через нее (а) не равна нулю; (б) положительна.

\smallskip

\noindent{\bf Задача \refstepcounter{problem}\arabic{problem}.
\label{odnorodnost}} Напряжение батарейки увеличили в $n$ раз. Докажите, что все силы тока в цепи также увеличились в $n$ раз.

%\smallskip
%\noindent{\bf Задача \refstepcounter{problem}\arabic{problem}.} В электрической цепи заменили батарейку (неизвестно, одинаково ли напряжение у старой и новой батарейки). Оказалось, что сила тока через батарейку не изменилась. Докажите, что силы тока через резисторы также не изменились.

\normalsize

\medskip

\textbf{Доказательство теоремы Дена о разрезании прямоугольника.}
\smallskip

%Мы готовы доказать теорему Дена.
Пусть прямоугольник разрезан на квадраты.
Расположим его так, чтобы две его стороны были вертикальны, а две другие --- горизонтальны. Будем считать, что длина горизонтальной стороны равна $1$. Ясно, что стороны всех квадратов либо вертикальны, либо горизонтальны. Рассмотрим электрическую цепь, соответствующую разрезанию.
Система уравнений, построенная по правилам Кирхгофа для этой цепи, имеет решение --- в качестве сил токов можно взять длины сторон квадратов. По теореме единственности
других решений у этой системы нет. Значит, по теореме о решении системы оно состоит из рациональных чисел. То есть длины сторон всех квадратов, а следовательно, и отношение сторон прямоугольника, рациональны.
%В частности, величина тока через батарейку рациональна.
%Эта величина численно равна вертикальной стороне прямоугольника.
Теорема Дена доказана.

\smallskip

\small

\noindent{\bf Задача \refstepcounter{problem}\arabic{problem}\label{root}.} Покажите, что квадрат нельзя разрезать на подобные (но не обязательно равные) прямоугольники с отношением сторон $\sqrt2$.
%\mscomm{Переместить после доказательства теоремы Дена?}
\normalsize

\smallskip

\bigskip
{\bf Десерт}
\smallskip

Мы ответили на все вопросы, поставленные в статье, но ее название осталось загадкой. Объяснение названия мы оставили на десерт: это будет наглядная картинка электрической цепи, построенной по разрезанию. Раньше физическая интерпретация выглядела, как некоторый трюк, теперь мы \emph{увидим, как до нее можно додуматься}.

%Представим себе прямоугольную металлическую пластинку, разбитую на квадратные. Вертикальные стороны большой пластинки соединим с полюсами батарейки. Вдоль горизонтальных разрезов изолируем квадратные пластинки  друг от друга, а вдоль вертикальных разрезов пусть они стыкуются. Каждая квадратная пластинка играет роль резистора. Известно, что сопротивление такого резистора не зависит от его размера; будем считать сопротивление равным $1$. Мы получили электрическую цепь. В действительности это та же цепь, что была у нас раньше.

%Наша физическая интерпретация может выглядеть как некоторый трюк. Мы оставили на десерт объяснение, как можно до нее додуматься.

%Как можно до нее додуматься, мы покажем в продолжении этой статьи в одном из следующих номеров журнала.

%\smallskip

%Выполним обещание, данное в конце предыдущей части статьи: покажем, как можно додуматься до физической интерпретации разрезания прямоугольника.

Пусть большой прямоугольник разрезан на меньшие (не обязательно квадраты), и требуется выразить
%Зададимся вопросом,
отношение сторон большого прямоугольника через отношения сторон меньших.
%Покажем на нескольких примерах, как найти отношение сторон прямоугольника, разрезанного на прямоугольники с заданными отношениями сторон.
% квадраты. Будем рассматривать разрезания не только на квадраты, но и на произвольные прямоугольники, отношения сторон которых считаются известными.
%\smallskip
%\small
%\refstepcounter{problem}
%\noindent{\bf Задача \arabic{problem}.} Большой прямоугольник разрезан на несколько меньших. Докажите, что стороны меньших прямоугольников параллельны сторонам большого. (\emph{Указание}. Если есть меньшие прямоугольники, расположенные ``криво'', то рассмотрите ближайший из них к вершине большого.)
%\normalsize
%\smallskip
%Очевидно, что в любом таком разрезании стороны всех рассматриваемых прямоугольников можно считать параллельными координатным осям, то есть либо {\it вертикальными}, либо {\it горизонтальными}. Под {\it отношением сторон} прямоугольника мы всегда понимаем отношение длины его горизонтальной стороны к длине вертикальной.
%%%%Будем обозначать отношение сторон прямоугольника так же, как и сам
%прямоугольник.
Расположим большой прямоугольник так, чтобы две его стороны были вертикальны, а две другие --- горизонтальны. %Очевидно, что тогда стороны всех меньших прямоугольников тоже либо вертикальны, либо горизонтальны.
{\it Отношением сторон} прямоугольника договоримся считать отношение длины его горизонтальной стороны к длине вертикальной.

%Но вернемся к обобщенной теореме Дена.
%Начнем с простейших примеров.

\smallskip

\begin{figure}[h]
\begin{tabular}{cc}
    \definecolor{zzttqq}{rgb}{0.6,0.2,0}
\begin{tikzpicture}[line cap=round,line join=round,>=triangle 45,x=0.5cm,y=0.5cm]
\clip(-3.78,0.84) rectangle (2.2,5.74);
\fill[color=zzttqq,fill=zzttqq,fill opacity=0.1] (-3,5) -- (2,5) -- (2,1) -- (-3,1) -- cycle;
\draw [color=zzttqq] (-3,5)-- (2,5);
\draw [color=zzttqq] (2,5)-- (2,1);
\draw [color=zzttqq] (2,1)-- (-3,1);
\draw [color=zzttqq] (-3,1)-- (-3,5);
\draw [color=zzttqq] (-1,5)-- (-1,1);
\draw[color=zzttqq] (0.6,3.2) node {$R_2$};
\draw[color=zzttqq] (-3.38,3.18) node {$x$};
\draw[color=zzttqq] (-1.86,3.18) node {$R_1$};
\end{tikzpicture}&\definecolor{zzttqq}{rgb}{0.6,0.2,0}
\begin{tikzpicture}[line cap=round,line join=round,>=triangle 45,x=0.5cm,y=0.5cm]
\clip(-3.78,0.76) rectangle (2.26,5.74);
\fill[color=zzttqq,fill=zzttqq,fill opacity=0.1] (-3,5) -- (2,5) -- (2,1) -- (-3,1) -- cycle;
\draw [color=zzttqq] (-3,5)-- (2,5);
\draw [color=zzttqq] (2,5)-- (2,1);
\draw [color=zzttqq] (2,1)-- (-3,1);
\draw [color=zzttqq] (-3,1)-- (-3,5);
\draw [color=zzttqq] (-3,3.34)-- (2,3.34);
\draw[color=zzttqq] (-0.74,5.36) node {$y$};
\draw[color=zzttqq] (-0.52,2.32) node {$R_2$};
\draw[color=zzttqq] (-0.56,4.26) node {$R_1$};
\end{tikzpicture}\\
    а & б
\end{tabular}
    \caption{Разрезания прямоугольника на $2$ прямоугольника}
    \label{ris2sq}
\end{figure}

\medskip

\noindent{\bf Пример 1.} Прямоугольник с отношением сторон
$R$ разделён вертикальным  разрезом на два прямоугольника
с отношениями сторон $R_1$ и $R_2$ (рис.~\ref{ris2sq}а). Покажем, что $R=R_1+R_2$. Действительно, пусть вертикальная сторона большого прямоугольника равна $x$. Тогда горизонтальные стороны меньших прямоугольников равны $R_1x$ и $R_2x$. Значит, $R=(R_1x+R_2x)/x=R_1+R_2$.

\smallskip

\noindent{\bf Пример 2.} Прямоугольник с отношением сторон
$R$ разделён горизонтальным разрезом на два прямоугольника
с отношениями сторон
$R_1$ и $R_2$ (рис.~\ref{ris2sq}б). Покажем, что $R=\frac{R_1R_2}{R_1+R_2}$.
Действительно, пусть горизонтальная сторона большого прямоугольника равна $y$.
Тогда вертикальные стороны меньших прямоугольников равны $y/R_1$ и
$y/R_2$.
Значит,
$$R=\frac{y}{y/R_1+y/R_2}=\frac{R_1R_2}{R_1+R_2}.$$
%Действительно, пусть вертикальные стороны меньших прямоугольников равны $I_1$ и
%$I_2$. Тогда горизонтальные стороны всех прямоугольников равны $I_1R_1=I_2R_2$.
%Значит,
%$$R=\frac{I_1R_1}{I_1+I_2}=\frac{R_1}{1+I_2/I_1}=\frac{R_1R_2}{R_2+R_2I_2/I_1}=\frac{R_1R_2}{R_1+R_2}.$$
%$R=\frac1{\frac1{R_1}+\frac1{R_2}}$. - такая форма записи не
%согласуется с теоремой о разрезании прямоугольника

%Покажем на примере, как искать отношение сторон прямоугольника в общей ситуации.

\smallskip

Да это же формулы сопротивления цепей из последовательно и параллельно соединённых резисторов (см.~рисунок~\ref{ris3})! Объяснение очень простое.
%, обратившись к теории электрических цепей;

\begin{figure}[h]
\begin{tabular}[t]{c}
$R=R_1+R_2$\\[3pt]
\definecolor{zzttqq}{rgb}{0.6,0.2,0}
\begin{tikzpicture}[line cap=round,line join=round,>=triangle 45,x=0.5cm,y=0.5cm]
\clip(-3.14,-2.22) rectangle (2.22,5.14);
\fill[color=zzttqq,fill=zzttqq,fill opacity=0.1] (-3,5) -- (2,5) -- (2,1) -- (-3,1) -- cycle;
\fill[color=zzttqq,fill=zzttqq,fill opacity=0.1] (-2.62,-2) -- (-2.62,-1.46) -- (-0.74,-1.46) -- (-0.74,-2) -- cycle;
\fill[color=zzttqq,fill=zzttqq,fill opacity=0.1] (1.62,-2) -- (1.62,-1.46) -- (-0.26,-1.46) -- (-0.26,-2) -- cycle;
\draw [color=zzttqq] (-3,5)-- (2,5);
\draw [color=zzttqq] (2,5)-- (2,1);
\draw [color=zzttqq] (2,1)-- (-3,1);
\draw [color=zzttqq] (-3,1)-- (-3,5);
\draw [color=zzttqq] (-1,5)-- (-1,1);
\draw [color=zzttqq] (-2.62,-2)-- (-2.62,-1.46);
\draw [color=zzttqq] (-2.62,-1.46)-- (-0.74,-1.46);
\draw [color=zzttqq] (-0.74,-1.46)-- (-0.74,-2);
\draw [color=zzttqq] (-0.74,-2)-- (-2.62,-2);
\draw [color=zzttqq] (1.62,-2)-- (1.62,-1.46);
\draw [color=zzttqq] (1.62,-1.46)-- (-0.26,-1.46);
\draw [color=zzttqq] (-0.26,-1.46)-- (-0.26,-2);
\draw [color=zzttqq] (-0.26,-2)-- (1.62,-2);
\draw (-3,-1.74)-- (-2.62,-1.74);
\draw [color=zzttqq] (-0.74,-1.72)-- (-0.26,-1.72);
\draw [color=zzttqq] (1.62,-1.72)-- (2,-1.72);
\draw [color=zzttqq] (-3,-1.74)-- (-3,-0.5);
\draw [color=zzttqq] (2,-1.72)-- (2,-0.5);
\draw [color=zzttqq] (-3,-0.5)-- (-0.54,-0.5);
\draw [color=zzttqq] (2,-0.5)-- (-0.26,-0.5);
\draw [color=zzttqq] (-0.54,-0.02)-- (-0.54,-1);
\draw [line width=2.4pt,color=zzttqq] (-0.26,-0.24)-- (-0.26,-0.82);
\draw[color=zzttqq] (0.6,3.2) node {$R_2$};
\draw[color=zzttqq] (-1.86,3.18) node {$R_1$};
\draw[color=zzttqq] (-1.7,-1.06) node {$R_1$};
\draw[color=zzttqq] (0.84,-1.02) node {$R_2$};
\end{tikzpicture}
\end{tabular}
\begin{tabular}[t]{c}
$R=\frac{R_1R_2}{R_1+R_2}$\\[2pt]
\definecolor{zzttqq}{rgb}{0.6,0.2,0}
\begin{tikzpicture}[line cap=round,line join=round,>=triangle 45,x=0.5cm,y=0.5cm]
\clip(-3.14,-4.02) rectangle (2.22,5.14);
\fill[color=zzttqq,fill=zzttqq,fill opacity=0.1] (-3,5) -- (2,5) -- (2,1) -- (-3,1) -- cycle;
\fill[color=zzttqq,fill=zzttqq,fill opacity=0.1] (-1.44,-2) -- (-1.44,-1.52) -- (0.52,-1.52) -- (0.52,-2) -- cycle;
\fill[color=zzttqq,fill=zzttqq,fill opacity=0.1] (-1.44,-3) -- (-1.44,-2.52) -- (0.52,-2.52) -- (0.52,-3) -- cycle;
\draw [color=zzttqq] (-3,5)-- (2,5);
\draw [color=zzttqq] (2,5)-- (2,1);
\draw [color=zzttqq] (2,1)-- (-3,1);
\draw [color=zzttqq] (-3,1)-- (-3,5);
\draw [color=zzttqq] (-1.44,-2)-- (-1.44,-1.52);
\draw [color=zzttqq] (-1.44,-1.52)-- (0.52,-1.52);
\draw [color=zzttqq] (0.52,-1.52)-- (0.52,-2);
\draw [color=zzttqq] (0.52,-2)-- (-1.44,-2);
\draw [color=zzttqq] (-3,-1.74)-- (-1.44,-1.74);
\draw [color=zzttqq] (-3,-1.74)-- (-3,-0.5);
\draw [color=zzttqq] (2,-1.72)-- (2,-0.5);
\draw [color=zzttqq] (-3,-0.5)-- (-0.54,-0.5);
\draw [color=zzttqq] (2,-0.5)-- (-0.26,-0.5);
\draw [color=zzttqq] (-0.54,-0.02)-- (-0.54,-1);
\draw [line width=2.4pt,color=zzttqq] (-0.26,-0.24)-- (-0.26,-0.82);
\draw [color=zzttqq] (-3,3)-- (2,3);
\draw [color=zzttqq] (2,-1.72)-- (0.52,-1.72);
\draw [color=zzttqq] (0.52,-3)-- (-1.44,-3);
\draw [color=zzttqq] (0.52,-2.52)-- (0.52,-3);
\draw [color=zzttqq] (-1.44,-2.52)-- (0.52,-2.52);
\draw [color=zzttqq] (-1.44,-3)-- (-1.44,-2.52);
\draw [color=zzttqq] (-3,-2.74)-- (-1.44,-2.74);
\draw [color=zzttqq] (2,-2.72)-- (0.52,-2.72);
\draw [color=zzttqq] (-1.44,-3)-- (-1.44,-2.52);
\draw [color=zzttqq] (-1.44,-2.52)-- (0.52,-2.52);
\draw [color=zzttqq] (0.52,-2.52)-- (0.52,-3);
\draw [color=zzttqq] (0.52,-3)-- (-1.44,-3);
\draw [color=zzttqq] (-3,-1.74)-- (-3,-2.74);
\draw [color=zzttqq] (2,-1.72)-- (2,-2.72);
\draw[color=zzttqq] (-0.36,2.14) node {$R_2$};
\draw[color=zzttqq] (-1.58,-1.02) node {$R_1$};
\draw[color=zzttqq] (-0.38,4.04) node {$R_1$};
\draw[color=zzttqq] (1.1,-3.52) node {$R_2$};
\end{tikzpicture}\\
\end{tabular}
%\begin{tabular}[t]{c}
%в\\
%\input{metal-fig5c.tex}
%\end{tabular}
\caption{
Формулы для отношения сторон такие же, как и для сопротивления!
}\label{ris3}
\end{figure}

\begin{figure}[h]
\definecolor{ffwwqq}{rgb}{1,0.4,0}
\definecolor{zzttqq}{rgb}{0.6,0.2,0}
\begin{tikzpicture}[line cap=round,line join=round,>=triangle 45,x=0.5cm,y=0.5cm]
\clip(-4.28,-1.3) rectangle (4.26,5.32);
\fill[color=zzttqq,fill=zzttqq,fill opacity=0.1] (-3,5) -- (-3,1) -- (-1,1) -- (-1,5) -- cycle;
\fill[color=zzttqq,fill=zzttqq,fill opacity=0.1] (0,5) -- (0,1) -- (3,1) -- (3,5) -- cycle;
\draw [color=zzttqq] (-3,-0.5)-- (-0.54,-0.5);
\draw [color=zzttqq] (2,-0.5)-- (-0.26,-0.5);
\draw [color=zzttqq] (-0.54,-0.02)-- (-0.54,-1);
\draw [line width=2.4pt,color=zzttqq] (-0.26,-0.24)-- (-0.26,-0.82);
\draw [line width=1.6pt,color=ffwwqq] (-3,5)-- (-3,1);
\draw [color=zzttqq] (-3,1)-- (-1,1);
\draw [line width=1.6pt,color=ffwwqq] (-1,1)-- (-1,5);
\draw [color=zzttqq] (-1,5)-- (-3,5);
\draw [line width=1.6pt,color=ffwwqq] (0,5)-- (0,1);
\draw [color=zzttqq] (0,1)-- (3,1);
\draw [line width=1.6pt,color=ffwwqq] (3,1)-- (3,5);
\draw [color=zzttqq] (3,5)-- (0,5);
\draw [color=zzttqq] (-2.98,-0.5)-- (-4,-0.5);
\draw [color=zzttqq] (-4,-0.5)-- (-4,3);
\draw (-4,3)-- (-3,3);
\draw [color=zzttqq] (-1,3)-- (0,3);
\draw [color=zzttqq] (3,3)-- (4,3);
\draw [color=zzttqq] (4,3)-- (4,-0.52);
\draw [color=zzttqq] (2,-0.5)-- (4,-0.52);
\draw[color=zzttqq] (-1.64,3.16) node {$R_1$};
\draw[color=zzttqq] (1.86,3.16) node {$R_2$};
\end{tikzpicture} \qquad
\definecolor{qqccqq}{rgb}{0,0.8,0}
\definecolor{ffwwqq}{rgb}{1,0.4,0}
\definecolor{zzttqq}{rgb}{0.6,0.2,0}
\begin{tikzpicture}[line cap=round,line join=round,>=triangle 45,x=0.5cm,y=0.5cm]
\clip(-4.28,-1.2) rectangle (3.28,6.32);
\fill[color=zzttqq,fill=zzttqq,fill opacity=0.1] (-3,1) -- (2,1) -- (2,3) -- (-3,3) -- cycle;
\fill[color=zzttqq,fill=zzttqq,fill opacity=0.1] (-3,6) -- (-3,4) -- (2,4) -- (2,6) -- cycle;
\draw [color=zzttqq] (-3,-0.5)-- (-0.54,-0.5);
\draw [color=zzttqq] (2,-0.5)-- (-0.26,-0.5);
\draw [color=zzttqq] (-0.54,-0.02)-- (-0.54,-1);
\draw [line width=2.4pt,color=zzttqq] (-0.26,-0.24)-- (-0.26,-0.82);
\draw [color=zzttqq] (-3,1)-- (2,1);
\draw [line width=1.6pt,color=ffwwqq] (2,1)-- (2,3);
\draw [line width=1.6pt,color=qqccqq] (2,3)-- (-3,3);
\draw [line width=1.6pt,color=ffwwqq] (-3,3)-- (-3,1);
\draw [line width=1.6pt,color=ffwwqq] (-3,6)-- (-3,4);
\draw [line width=1.6pt,color=qqccqq] (-3,4)-- (2,4);
\draw [line width=1.6pt,color=ffwwqq] (2,4)-- (2,6);
\draw [color=zzttqq] (2,6)-- (-3,6);
\draw [color=zzttqq] (-3,-0.5)-- (-4,-0.5);
\draw [color=zzttqq] (-4,-0.5)-- (-4,5);
\draw [color=zzttqq] (-4,5)-- (-3,5);
\draw [color=zzttqq] (-3,2)-- (-4,2);
\draw [color=zzttqq] (2,-0.5)-- (3,-0.5);
\draw [color=zzttqq] (3,-0.5)-- (3,5);
\draw [color=zzttqq] (3,5)-- (2,5);
\draw [color=zzttqq] (2,2)-- (3,2);
\draw[color=zzttqq] (-0.14,2.16) node {$R_2$};
\draw[color=zzttqq] (-0.14,5.16) node {$R_1$};
\end{tikzpicture} \qquad
\definecolor{uququq}{rgb}{0.25,0.25,0.25}
\definecolor{zzttqq}{rgb}{0.6,0.2,0}
\begin{tikzpicture}[line cap=round,line join=round,>=triangle 45,x=0.5cm,y=0.5cm]
\clip(-10.38,-19.28) rectangle (10.26,1.33);
\fill[color=zzttqq,fill=zzttqq,fill opacity=0.1] (-9,-0.5) -- (-4.5,-0.5) -- (-4.5,-5) -- (-9,-5) -- cycle;
\fill[color=zzttqq,fill=zzttqq,fill opacity=0.1] (-9,-6) -- (-4,-6) -- (-4,-11) -- (-9,-11) -- cycle;
\fill[color=zzttqq,fill=zzttqq,fill opacity=0.1] (-3.5,-0.5) -- (-3.5,-4.5) -- (0.5,-4.5) -- (0.5,-0.5) -- cycle;
\fill[color=zzttqq,fill=zzttqq,fill opacity=0.1] (-4,-5) -- (-3.5,-5) -- (-3.5,-5.5) -- (-4,-5.5) -- cycle;
\fill[color=zzttqq,fill=zzttqq,fill opacity=0.1] (-3,-5.5) -- (0.5,-5.5) -- (0.5,-9) -- (-3,-9) -- cycle;
\fill[color=zzttqq,fill=zzttqq,fill opacity=0.1] (-3,-9.5) -- (-1,-9.5) -- (-1,-11.5) -- (-3,-11.5) -- cycle;
\fill[color=zzttqq,fill=zzttqq,fill opacity=0.1] (-9,-12) -- (-2,-12) -- (-2,-19) -- (-9,-19) -- cycle;
\fill[color=zzttqq,fill=zzttqq,fill opacity=0.1] (1.5,-0.5) -- (1.5,-8) -- (9,-8) -- (9,-0.5) -- cycle;
\fill[color=zzttqq,fill=zzttqq,fill opacity=0.1] (9,-19) -- (0,-19) -- (0,-10) -- (9,-10) -- cycle;
\draw [color=zzttqq] (-2,0.56)-- (0.32,0.56);
\draw [color=zzttqq] (3,0.56)-- (0.52,0.56);
\draw [color=zzttqq] (0.32,1.06)-- (0.32,0.06);
\draw [line width=2.4pt,color=zzttqq] (0.52,0.8)-- (0.52,0.32);
\draw [color=zzttqq] (-9,-0.5)-- (-4.5,-0.5);
\draw [color=zzttqq] (-4.5,-0.5)-- (-4.5,-5);
\draw [color=zzttqq] (-4.5,-5)-- (-9,-5);
\draw [color=zzttqq] (-9,-5)-- (-9,-0.5);
\draw [color=zzttqq] (-9,-6)-- (-4,-6);
\draw [color=zzttqq] (-4,-6)-- (-4,-11);
\draw [color=zzttqq] (-4,-11)-- (-9,-11);
\draw [color=zzttqq] (-9,-11)-- (-9,-6);
\draw [color=zzttqq] (-3.5,-0.5)-- (-3.5,-4.5);
\draw [color=zzttqq] (-3.5,-4.5)-- (0.5,-4.5);
\draw [color=zzttqq] (0.5,-4.5)-- (0.5,-0.5);
\draw [color=zzttqq] (0.5,-0.5)-- (-3.5,-0.5);
\draw [color=zzttqq] (-4,-5)-- (-3.5,-5);
\draw [color=zzttqq] (-3.5,-5)-- (-3.5,-5.5);
\draw [color=zzttqq] (-3.5,-5.5)-- (-4,-5.5);
\draw [color=zzttqq] (-4,-5.5)-- (-4,-5);
\draw [color=zzttqq] (-3,-5.5)-- (0.5,-5.5);
\draw [color=zzttqq] (0.5,-5.5)-- (0.5,-9);
\draw [color=zzttqq] (0.5,-9)-- (-3,-9);
\draw [color=zzttqq] (-3,-9)-- (-3,-5.5);
\draw [color=zzttqq] (-3,-9.5)-- (-1,-9.5);
\draw [color=zzttqq] (-1,-9.5)-- (-1,-11.5);
\draw [color=zzttqq] (-1,-11.5)-- (-3,-11.5);
\draw [color=zzttqq] (-3,-11.5)-- (-3,-9.5);
\draw [color=zzttqq] (-9,-12)-- (-2,-12);
\draw [color=zzttqq] (-2,-12)-- (-2,-19);
\draw [color=zzttqq] (-2,-19)-- (-9,-19);
\draw [color=zzttqq] (-9,-19)-- (-9,-12);
\draw [color=zzttqq] (1.5,-0.5)-- (1.5,-8);
\draw [color=zzttqq] (1.5,-8)-- (9,-8);
\draw [color=zzttqq] (9,-8)-- (9,-0.5);
\draw [color=zzttqq] (9,-0.5)-- (1.5,-0.5);
\draw [color=zzttqq] (9,-19)-- (0,-19);
\draw [color=zzttqq] (0,-19)-- (0,-10);
\draw [color=zzttqq] (0,-10)-- (9,-10);
\draw [color=zzttqq] (9,-10)-- (9,-19);
\draw [color=zzttqq] (-2,0.56)-- (-10.02,0.58);
\draw [color=zzttqq] (3,0.56)-- (9.99,0.58);
\draw [color=zzttqq] (-10.02,0.58)-- (-10,-15.5);
\draw (-10,-15.5)-- (-9,-15.5);
\draw [color=zzttqq] (-10.01,-8.5)-- (-9,-8.5);
\draw [color=zzttqq] (-10.01,-2.72)-- (-9,-2.72);
\draw [color=zzttqq] (-4.5,-2.5)-- (-3.5,-2.5);
\draw [color=zzttqq] (9.99,0.58)-- (10,-14.5);
\draw [color=zzttqq] (10,-14.5)-- (9,-14.5);
\draw [color=zzttqq] (-2,-15)-- (0,-15);
\draw [color=zzttqq] (9,-4.19)-- (10,-4.21);
\draw [color=zzttqq] (1.5,-2.5)-- (1,-2.5);
\draw [color=zzttqq] (0.5,-2.5)-- (1,-2.5);
\draw [color=zzttqq] (1,-2.5)-- (0.99,-7.16);
\draw [color=zzttqq] (0.99,-7.16)-- (0.5,-7.19);
\draw [color=zzttqq] (-3,-7.5)-- (-4,-7.5);
\draw [color=zzttqq] (-3.25,-5.24)-- (-3.27,-10.49);
\draw [color=zzttqq] (-3.27,-10.49)-- (-3,-10.5);
\draw [color=zzttqq] (-3.25,-5.24)-- (-3.5,-5.23);
\draw [color=zzttqq] (-4.27,-2.5)-- (-4.27,-5.25);
\draw [color=zzttqq] (-4.27,-5.25)-- (-4,-5.25);
\draw [color=zzttqq] (-1,-10.5)-- (-0.5,-10.5);
\draw [color=zzttqq] (-0.5,-10.5)-- (-0.51,-15);
\fill [color=uququq] (-3.26,-7.5) ++(-1.0pt,0 pt) -- ++(1.0pt,1.0pt)--++(1.0pt,-1.0pt)--++(-1.0pt,-1.0pt)--++(-1.0pt,1.0pt);
\end{tikzpicture}
\caption{Электрические цепи из металлических пластинок}
\label{plates}
\end{figure}

Представим себе, что у нас есть прямоугольная металлическая пластинка. Соединим ее вертикальные стороны с полюсами батарейки (точнее, к каждой из вертикальных сторон по всей длине приложим проводник, соединенный с соответствующим полюсом). Тогда через пластинку пойдет ток в горизонтальном направлении. Пластинка играет роль резистора. Как известно из физики, ее сопротивление пропорционально отношению длины к площади вертикального поперечного сечения. Иными словами, сопротивление пластинки пропорционально отношению ее сторон.
%Выбрав подходящую систему единиц измерения, можно
Для простоты будем считать коэффициент пропорциональности равным $1$.

%``Физическое'' объяснение формулы из примера 1 такое.
Для примера 1 приставим друг к другу две прямоугольные пластинки одинаковыми вертикальными сторонами, см. рисунок~\ref{plates} слева вверху. Оставшиеся вертикальные стороны соединим с полюсами батарейки. Получим цепь из двух последовательно соединенных резисторов. Если отношения сторон этих пластинок --- $R_1$ и $R_2$, то их сопротивления  --- тоже $R_1$ и $R_2$. Сопротивление большой пластинки, составленной из двух, равно $R_1+R_2$, как сопротивление двух последовательно соединенных резисторов. Вот ``физическое'' объяснение того, что отношение сторон большой пластинки равно $R_1+R_2$.

Перейдем к примеру 2. Приставим две пластинки друг к другу одинаковыми горизонтальными сторонами, а вертикальные стороны соединим с полюсами батарейки,
см. рисунок~\ref{plates} справа вверху.
Поскольку ток течет в горизонтальном направлении, то через линию стыковки ток не идет. Изолируем пластинки друг от друга: ток и сопротивление цепи не изменятся. Мы получим пару параллельно соединенных резисторов, значит, отношение сторон большого прямоугольника находится по формуле $\frac{R_1R_2}{R_1+R_2}$.

Этот же метод можно применить для любого разрезания прямоугольника на прямоугольники, скажем, для изображенного на рисунках~\ref{ris111} и~\ref{plates} снизу. Представим себе большую прямоугольную пластинку, разбитую на меньшие. Вертикальные стороны большой пластинки соединим с полюсами батарейки. Вдоль горизонтальных линий разреза изолируем меньшие пластинки  друг от друга, а вдоль вертикальных линий разреза пусть они стыкуются. %Каждая меньшая пластинка будет играть роль резистора. %Как известно из физики, сопротивление маленькой пластинки пропорционально отношению длины к площади поперечного сечения, то есть, пропорционально отношению сторон маленькой пластинки.
%Общее сопротивление полученной электрической цепи равно сопротивлению большой пластинки, то есть, пропорциональному отношению сторон большой пластинки.
%Так же, как и раньше, отношение сторон большой пластинки можно выразить через отношения сторон малых пластинок, пользуясь некоторой формулой для общего сопротивления электрической цепи.
Мы получили электрическую цепь, сопротивление которой равно отношению сторон большой пластинки.
На самом деле это та же самая электрическая цепь, которую мы построили раньше.

Таким образом, чтобы найти отношение сторон прямоугольника, достаточно измерить сопротивление электрической цепи, построенной по разрезанию! Об этом мы подробно расскажем в одном из следующих номеров журнала.

%\mscomm{Дать определение сопротивления цепи и несколько задач про него?}

\smallskip

\small
\noindent{\bf Задача \refstepcounter{problem}\arabic{problem}\label{five}.} Прямоугольник разделён на пять прямоугольников с отношениями сторон $R_1=R_2=R_3=1$, $R_4=R_5=3$ так, как показано на рисунке~\ref{ris2222}. Найдите отношение сторон большого прямоугольника.

\normalsize

\begin{figure}[h]
    \definecolor{zzttqq}{rgb}{0.6,0.2,0}
\begin{tikzpicture}[line cap=round,line join=round,>=triangle 45,x=0.7cm,y=0.7cm]
\clip(-2.24,1.76) rectangle (3.22,5.16);
\fill[color=zzttqq,fill=zzttqq,fill opacity=0.1] (-2,5) -- (3,5) -- (3,2) -- (-2,2) -- cycle;
\draw [color=zzttqq] (-2,5)-- (3,5);
\draw [color=zzttqq] (3,5)-- (3,2);
\draw [color=zzttqq] (3,2)-- (-2,2);
\draw [color=zzttqq] (-2,2)-- (-2,5);
\draw [color=zzttqq] (0,5)-- (0,3);
\draw [color=zzttqq] (-2,3)-- (1,3);
\draw [color=zzttqq] (1,4)-- (1,2);
\draw [color=zzttqq] (0,4)-- (3,4);
\draw[color=zzttqq] (1.68,4.6) node {$R_4$};
\draw[color=zzttqq] (2.16,3.16) node {$R_3$};
\draw[color=zzttqq] (-0.36,2.54) node {$R_5$};
\draw[color=zzttqq] (-1.02,4.08) node {$R_1$};
\draw[color=zzttqq] (0.58,3.6) node {$R_2$};
\end{tikzpicture}
    \caption{Разрезание прямоугольника на $5$ прямоугольников} \label{ris2222}
\end{figure}

\small
\smallskip

\noindent{\bf Задача \refstepcounter{problem}\arabic{problem}$^{**}$\label{chain2rect}.}
%Пусть электрическая цепь нарисована на плоскости так, что ее ребра не имеют общих внутренних точек.
%Предположим, что  в любой грани нет клемм с одинаковым потенциалом.
На плоскости дана электрическая цепь из резисторов сопротивлением $1$ и батарейки напряжением $1$. Предположим, все силы токов в цепи ненулевые, а все контуры --- треугольные\footnote{На самом деле, достаточно предположения, что в каждом контуре \emph{потенциалы} вершин попарно различны.}.
Попробуйте доказать, что тогда эта цепь получается из некоторого разрезания прямоугольника на квадраты.

\normalsize

\medskip

{\bf Благодарности.} Авторы благодарны
Евгению Выродову, Сергею Маркелову, Евгению Могилевскому, Владимиру Протасову, Александру Прохорову, Святославу Фельдшерову и Борису Френкину за ценные замечания.

\bigskip
\bigskip

\textbf{Решения задач.}

\textbf{\ref{plane}.} Решим сразу пункт б). Возьмем пример на рисунке 3 статьи, где стороны всех квадратов целые: например, можно взять разбиение, в котором прямоугольник, разбитый на квадраты, имеет размеры $33\times 32$ (первой указывается длина вертикальной стороны, потом горизонтальной). Будем замощать плоскость далее по такому алгоритму. Приставляем к нашему прямоугольнику слева квадрат $33\times 33$ (сторона к стороне), к получившемуся прямоугольнику размером $33\times 65$ приставляем сверху квадрат $65\times 65$ (сторона к стороне), к получившемуся прямоугольнику $98\times 65$ приставляем справа квадрат $98\times 98$, затем приставляем снизу квадрат $163\times 163$, снова приставляем квадрат слева, и так далее. В результате мы заполняем плоскость «по спирали», и каждый новый квадрат имеет сторону, равную большей стороне прямоугольника, к которому он приставляется – тем самым он больше по размеру, чем все предыдущие квадраты.

См. также F.V. Henle, J.M. Henle, \textrm{Squaring the plane}, Amer.~Math.~Monthly \textrm{115:1} (2008), 3--12; \url{http://citeseerx.ist.psu.edu/viewdoc/download?doi=10.1.1.138.7633&rep=rep1&type=pdf}.

\textbf{\ref{cube}.} Нельзя.

\emph{Указание}. Рассмотрите кубик наименьшего размера, примыкающий к одной из граней большого куба. Покажите, что к нему примыкает кубик еще меньшего размера, к тому – еще меньшего, и т.д.
Подробное решение можно найти в статье Л. Курляндчика и Г. Розенблюма ``Метод бесконечного спуска'' в «Кванте» \No 1 за 1978 год.

\textbf{\ref{rooms}.}
Нет.

Пусть это возможно. Из плана квартиры (рисунок~\ref{ris9} статьи) видно, что сторона левой верхней комнаты длиннее, чем сторона правой верхней, та, в свою очередь, длиннее стороны правой нижней комнаты, та длиннее стороны левой нижней, и, наконец, последняя длиннее стороны левой верхней. Получили, что сторона левой верхней комнаты длиннее самой себя, что невозможно.

\textbf{\ref{closed chain}.} Рассмотрим замкнутую цепочку ребер, не проходящую ни через какую вершину дважды,
см.~рисунок~\ref{cepochka}. Она ограничивает некоторый многоугольник на плоскости. Запишем второе правило Кирхгофа для всех контуров, которые попали внутрь многоугольника, и сложим полученные уравнения. Тогда напряжения для всех ребер, не лежащих на границе многоугольника, сократятся. Действительно, каждое такое ребро входит в два контура: для одного из них направление стрелки на ребре совпадает с направлением обхода, для другого --- противоположно направлению обхода. Для первого контура напряжение на ребре будет входить в уравнение со знаком плюс, а для второго --- со знаком минус. Значит, сумма всех выписанных уравнений --- это в точности второе правило Кирхгофа, только роль контура играет наша замкнутая цепочка.

\begin{figure}[htb]
\definecolor{ffqqqq}{rgb}{1,0,0}
\definecolor{qqwwff}{rgb}{0,0.4,1}
\begin{tikzpicture}[line cap=round,line join=round,>=triangle 45,x=0.5cm,y=0.5cm]
\clip(-4.12,-1.74) rectangle (8.36,5.66);
\draw [color=qqwwff] (6.9,1.08)-- (1.5,-1.2);
\draw [color=qqwwff] (6.9,1.08)-- (2.4,3.76);
\draw [color=qqwwff] (1.5,-1.2)-- (-2.7,1.04);
\draw [color=qqwwff] (0.3,1.8)-- (2.4,3.76);
\draw [color=qqwwff] (1.5,-1.2)-- (0.3,1.8);
\draw [color=qqwwff] (0,4.66)-- (2.4,3.76);
\draw [color=qqwwff] (0,4.66)-- (-2.7,1.04);
\draw [color=qqwwff] (0.3,1.8)-- (-2.7,1.04);
\draw [color=qqwwff] (0,4.66)-- (0.3,1.8);
\draw [->,color=qqwwff] (-2.7,1.04) -- (0,4.66);
\draw [->,color=qqwwff] (0,4.66) -- (2.4,3.76);
\draw [->,color=qqwwff] (2.4,3.76) -- (6.9,1.08);
\draw [->,color=qqwwff] (0,4.66) -- (0.3,1.8);
\draw [->,color=qqwwff] (-2.7,1.04) -- (0.3,1.8);
\draw [->,color=qqwwff] (0.3,1.8) -- (2.4,3.76);
\draw [->,color=qqwwff] (-2.7,1.04) -- (1.5,-1.2);
\draw [->,color=qqwwff] (0.3,1.8) -- (1.5,-1.2);
\draw [->,color=qqwwff] (1.5,-1.2) -- (6.9,1.08);
\draw [dash pattern=on 2pt off 2pt] (2.92,1.42) circle (0.59cm);
\draw [dash pattern=on 2pt off 2pt] (0.88,3.48) circle (0.18cm);
\draw [dash pattern=on 2pt off 2pt] (-0.68,2.62) circle (0.22cm);
\draw [dash pattern=on 2pt off 2pt] (-0.28,0.7) circle (0.27cm);
\draw [->] (0.9,3.84) -- (1.15,3.72);
\draw [->] (3.06,2.59) -- (3.46,2.47);
\draw [->] (-0.3,1.24) -- (0.02,1.15);
\draw [->] (-0.68,3.06) -- (-0.43,2.99);
\draw[color=qqwwff] (4.28,-0.6) node {7};
\fill [color=ffqqqq] (6.9,1.08) circle (2.0pt);
\draw[color=ffqqqq] (7.34,1.18) node {$-$};
\draw[color=qqwwff] (5,2.98) node {9};
\draw[color=qqwwff] (-0.88,-0.86) node {8};
\draw[color=qqwwff] (1.7,2.48) node {5};
\fill [color=ffqqqq] (1.5,-1.2) circle (2.0pt);
\draw[color=qqwwff] (1.28,0.8) node {6};
\fill [color=ffqqqq] (2.4,3.76) circle (2.0pt);
\draw[color=qqwwff] (1.38,4.84) node {4};
\draw[color=qqwwff] (-1.32,4.04) node {3};
\fill [color=ffqqqq] (-2.7,1.04) circle (2.0pt);
\draw[color=ffqqqq] (-3.22,1.08) node {$+$};
\draw[color=qqwwff] (-1.36,1.84) node {2};
\fill [color=ffqqqq] (0,4.66) circle (2.0pt);
\fill [color=ffqqqq] (0.3,1.8) circle (2.0pt);
\draw[color=qqwwff] (-0.26,3.78) node {1};
\end{tikzpicture}
\caption{Второе правило Кирхгофа для замкнутой цепочки ребер --- это сумма вторых правил Кирхгофа для всех контуров внутри нее.}
\label{cepochka}
\end{figure}

\textbf{\ref{polozhitelnost}.} Докажем сразу пункт б). Сначала заметим, что в цепи есть ненулевые токи. Действительно, запишем второе правило Кирхгофа для контура, содержащего батарейку. В правой части уравнения стоит напряжение батарейки, ненулевое число. Значит, одно из слагаемых в левой части не равно нулю, то есть, сила тока через одно из ребер ненулевая.

Теперь будем рассуждать, как в доказательстве принципа техники безопасности. Поменяем одновременно направление стрелки, знак силы тока и знак напряжения на каждом ребре с отрицательной силой тока. В частности, если ток через батарейку был отрицателен, то после нашей замены напряжение на ней станет отрицательно. Правила Кирхгофа по-прежнему будут выполняться. Начнем движение с ребра, на котором сила тока ненулевая, и будем двигаться в направлении стрелок, не заходя в ребра с нулевой силой тока. Мы можем неограниченно продолжать движение, и рано или поздно мы впервые получим замкнутую цепочку ребер.

Запишем второе правило Кирхгофа для этой цепочки, см.~задачу~\ref{closed chain}. В левой части стоит положительное число. Значит, в правой части тоже должно стоять положительное число. Но это возможно, только если контур проходит через батарейку. Поскольку мы двигались по ребрам с положительной силой тока, то ток через батарейку положителен, в частности, ненулевой. Как мы отметили выше, если бы мы меняли направление стрелки на батарейке, то напряжение на ней стало бы отрицательно, то есть в правой части выписанного уравнения стояло бы отрицательное число. Значит, ток через батарейку с самого начала (до смены направлений стрелок) был положителен.

\textbf{\ref{odnorodnost}.} Пусть силы тока в цепи до замены батарейки были $I_1,I_2,\dots$. Тогда $I_1,I_2,\dots$ --- решение системы уравнений, построенной по правилам Кирхгофа. Умножим каждое уравнение на число $n$. Получим, что $nI_1,nI_2,\dots$ --- решение системы уравнений, полученной увеличением напряжения батарейки в $n$ раз. Новая система соответствует правилам Кирхгофа после замены батарейки. По теореме единственности новая система других решений не имеет. Значит, $nI_1,nI_2,\dots$ --- и есть силы тока в новой цепи.

\textbf{\ref{root}.} Предположим, что такое разрезание возможно. Растянем наш квадрат в $\sqrt{2}$ раз вдоль одной из его сторон, он превратится в прямоугольник с отношением сторон $\sqrt{2}$. А прямоугольники, на которые он был разрезан, превратятся либо в квадраты, либо в прямоугольники с отношением сторон $2$, каждый из которых можно разрезать на два квадрата! Получили, что прямоугольник с иррациональным отношением сторон (наш растянутый квадрат) разрезается на квадраты, а это противоречит теореме Дена.

\textbf{\ref{five}.} $5:3$

Конечно же, можно выписать уравнения стыковки и решить систему. Но можно просто угадать ответ: видно, что если взять прямоугольники $4$ и $5$ размером $1\times 3$, прямоугольники $1$ и $3$ взять размером $2\times 2$, а прямоугольник $2$ взять размером $1\times 1$, то мы получим разрезание как на картинке, и искомое отношение сторон будет равно $5:3$. По теореме единственности это и есть ответ.

\textbf{\ref{chain2rect}.} Известное нам решение этой задачи очень сложное, и мы его не приводим.

\end{document}